\documentclass[aos,MSNbibl,nameyear,dvips]{arximspdf}
\usepackage{graphicx}

\doi{10.1214/12-AOS976} 
\volume{40}
\issue{3}
\pubyear{2012}
\firstpage{1430}
\lastpage{1464}

\makeatletter

\newtheorem{theorem}{Theorem}

\newproclaim{definition}{Definition}
\newproclaim{example}{Example}

\newtheorem{lemma}{Lemma}
\newtheorem{proposition}{Proposition}

\newproclaim{remark}{Remark}
\newproclaim{assumption}{Assumption}
\newproclaim{Extension}{Extension}
\newproclaim{Step}{Step}

\makeatother

\begin{document}
\begin{frontmatter}

\title{Identifying the successive Blumenthal--Getoor indices of
a discretely observed process}
\runtitle{Successive Blumenthal--Getoor indices}

\begin{aug}
\author[A]{\fnms{Yacine} \snm{A\"{i}t-Sahalia}\corref{}\thanksref{t1}\ead[label=e1]{yacine@princeton.edu}}
\and
\author[B]{\fnms{Jean} \snm{Jacod}\ead[label=e2]{jean.jacod@upmc.fr}}
\runauthor{Y. A\"{i}t-Sahalia and J. Jacod}
\affiliation{Princeton University and Universit\`e Pierre and Marie Curie}
\address[A]{Department of Economics\\
Princeton University and NBER\\
Princeton, New Jersey 08544-1021\\
USA\\
\printead{e1}}
\address[B]{Institut de Math\'ematiques de Jussieu\\
CNRS UMR 7586\\
Universit\`e Pierre and Marie Curie\\
75252 Paris C\'edex 05\\
France\\
\printead{e2}} 
\end{aug}

\thankstext{t1}{Supported in part by NSF Grant SES-0850533.}

\received{\smonth{4} \syear{2011}}
\revised{\smonth{1} \syear{2012}}

%
\begin{abstract}
This paper studies the identification of the L\'{e}vy jump measure of a
discretely-sampled semimartingale. We define successive
Blumenthal--Getoor indices of jump activity, and show that the leading
index can always be identified, but that higher order indices are only
identifiable if they are sufficiently close to the previous one, even
if the path is fully observed. This result establishes a clear boundary
on which aspects of the jump measure can be identified on the basis of
discrete observations, and which cannot. We then propose an estimation
procedure for the identifiable indices and compare the rates of
convergence of these estimators with the optimal rates in a special
parametric case, which we can compute explicitly.
\end{abstract}

%

\end{frontmatter}

\section{Introduction}\label{secintro}

Let $X$ be a one-dimensional semimartingale defined on a finite time interval
$[0,T]$. Our objective is to make some progress toward the
identification of
the jump measure of $X$ at high frequency. The motivation for what
follows has
its roots in a family of econometric problems, which can be stated as follows.
We observe a single path of $X$, but not fully: although other observation
schemes are possible, the most typical is one where we observe the variables
$X_{i\Delta_{n}}$ for $i=0,1,\ldots,[T/\Delta_{n}]$, where $[x]$
denotes the
integer part of the real $x$, over a fixed observation span $T$ and where
$\Delta_{n}$ is small. Asymptotic results are derived in the high-frequency
limit where the sequence $\Delta_{n}$ going to $0$. The overall
objective is
to find out what can be recovered, that is, identified, about the
dynamics of
$X$, in this setup where a single path, partially observed at a
discrete time
interval, is all that is available. For those parameters which can be
identified, we also want asymptotically consistent estimators, with a rate
whenever possible.

For the dynamics of $X$, we restrict our attention to It\^{o} semimartingales,
meaning that the characteristics $(B,C,\nu)$ of $X$ can be written as
follows:
%
\begin{eqnarray} \label{3}%
B_{t}(\omega)&=&\int_{0}^{t}b_{s}( \omega) \,ds,\nonumber\\[-8pt]\\[-8pt]
C_{t}(\omega)&=&\int_{0}^{t}c_{s}(\omega)\,ds,\qquad \nu(\omega,dt,dx)=dt\otimes
F_{t}(\omega,dx)\nonumber
\end{eqnarray}
for some adapted processes $b_{t}$ and $c_{t}$ and measure $F_{t}(\omega,dx)$.
Recall that~$B$ is the drift, $C$ is the quadratic variation of the continuous
martingale part and~$\nu$ is the compensator of the jump measure $\mu$
of $X$
[see \citet{jacodshiryaev2003} for more details on characteristics]. As
is well
known, these are the canonical models for arbitrage-free asset prices.

A sizeable part of the paper, however, is concerned with the much-restrict\-ed
class of L\'evy processes. A semimartingale $X$ is a L\'evy process if and
only if~(\ref{3}) holds with $b_{t}(\omega)=b\in\mathbb{R}$ and $c_{t}%
(\omega)=c\geq0$ and $F_{t}(\omega,dx)=F(dx)$ independent of $\omega$
and $t$.
The measure $F$ is the \textit{L\'evy measure}, and it integrates
$x^{2}\wedge
1$. The (deterministic) triple $(b,c,F)$ is then the characteristic triple
coming in the L\'evy--Khintchine formula, providing the characteristic function
of $X_{t}$,
%
\begin{equation} \label{1}
\mathbb{E}[ e^{iuX_{t}}] =\exp t\biggl( iub-\frac{cu^{2}}{2}%
+\int\bigl( e^{iux}-1-iux1_{\{|x|\leq1\}}\bigr) F(dx)\biggr).
\end{equation}
This completely characterizes the entire law of $X$.

Ultimately, we would like to identify as much as we can of the characteristics
$B$, $C$ and $\nu$, and give consistent estimators for the identifiable
parameters. The situation is well understood for the first two
characteristics, $B$ and $C$. When $X$ is fully observed on $[0,T]$,
one knows
the jumps (size and location) occurring within the interval, and the quadratic
variation of $X$ on $[0,T]$, hence the function $t\mapsto C_{t}$ on $[0,T]$.
On the other hand, and at least when $C$ is strictly increasing (which
is the
case in almost all models used in practice), nothing can be said about the
drift~$B$. When the process is observed only at discrete times, $C_{t}$
is no
longer exactly known, but there are well established methods to
estimate it in
a consistent way as the observation mesh goes to $0$, even in the
presence of jumps.

We focus on the remaining open question, which concerns identifiability and
estimation for the third characteristic, $\nu$, or equivalently
$F_{t}$, for a
discretely sampled semimartingale. The measure $F_{t}$ in a sense describes
the law of a jump occurring at time $t$, conditionally on the past before
$t$. There is a~vast literature on identifying the L\'{e}vy measure
when the
time horizon $T$ is asymptotically infinite, and when $X$ is a L\'evy process;
see, for example, \citet{basawabrockwell82}, \citet{figueroalopezhoudre06},
\citet{nishiyama09}, \citet{neumannreiss09} and \citet
{comtegenoncatalot09}. But
over a finite time horizon $T$, we cannot reconstruct $\nu$ fully because
there are only finitely many jumps on $[0,T]$ with size bigger than any
$\varepsilon>0$. The open question which we seek to address in this
paper is:
what can we and can we not identify about $\nu?$ High-frequency data analysis
has proved a very fruitful area of research. As we will see, however,
it is
not able to achieve everything, and our objective in this paper is to pinpoint
exactly the limitations, or frontier, involved in using high-frequency data
over a fixed time span.

We can say something about the concentration of $\nu$ around $0$. For example,
we can decide for which $p\geq0$ we have $\int_{0}^{T}ds\int F_{t}%
(\omega,dx)(|x|^{p}\wedge1)<\infty$, because outside a null set again these
are exactly those $p$'s for which $\sum_{s\leq T}|\Delta X_{s}(\omega
)|^{p}<\infty$,\vspace*{1pt} where $\Delta X_{s}=X_{s}-X_{s-}$ is the size of the
jump at
time~$s$, if any. The infimum of all such $p$'s is a generalization of the
\textit{Blumenthal--Getoor} index (or BG index) of the process up to time $T$,
and it is \textit{known} when $X$ is fully observed. Note that a
priori it is random, and also increasing with $T$, and always with
values in
$[0,2]$. However, in the L\'{e}vy process case, it reduces to $\inf\{
p\dvtx\int
F(dx)(|x|^{p}\wedge1)<\infty\}$, and is nonrandom and independent of
time. It
was originally introduced by \citet{blumenthalgetoor61}, and for a stable
process the BG index is also the stability index of the process.

The interest in identifying the BG index lies in the fact that the index
allows for a classification of the processes from least active to most active:
processes with BG index equal to $0$ are either finitely active or infinitely
active but with slow, sub-polynomial, divergence of $\nu$ near $0;$ processes
with BG index strictly positive are all infinitely active; processes
with BG
index less than $1$ have paths of finite variation; processes with BG index
greater than $1$ have paths of infinite variation; and in the limit, processes
with continuous paths have an ``activity index'' (the analog of the BG index
which no longer exists) equal to $2$ when the volatility is not
vanishing. In
other words, jumps become more and more active as the BG index
increases from
$0$ to $2$, and we can think of this generalized BG index as an \textit{index
of jump activity}.

In the case of discrete observations at times $i\Delta_{n}$ with $\Delta_{n}$
going to $0$, recovering the random BG index in full generality seems
out of
reach, but \citet{yacjacod09b} constructed estimators of the nonrandom number
$\beta$ that are consistent as $\Delta_{n}\rightarrow0$, under the main
assumption that locally near $0$, we have the behavior
%
\begin{equation} \label{I1}%
F_{t}(\omega,[-u,u]^{c})\sim\frac{a_{t}(\omega)}{u^{\beta}}\qquad
\mbox{as }u\downarrow0
\end{equation}
(plus a few technical hypotheses), where $a_{t}\geq0$ is a process: in this
case, $\beta$ is the---deterministic---BG index at time $t$, on the set
$\{ \int_{0}^{t}a_{s}\,ds>0\} $. We call this behavior
``proto-stable,'' since it is similar
to that
of a stable process but only near~$0$. Away from a neighborhood of $0$, the
jump measure is completely unrestricted. We obtained the rate of convergence
and a central limit theorem for the estimators, depending upon the rate
in the
approximations~(\ref{I1}). Related estimators or tests for $\beta$ include
\citet{belomestny10}, \citet{contmancini11} and \citet{tauchentodorov10}.

We can think of~(\ref{I1}) as providing the leading term, near $0$, of the
jump measure of $X$. Given that this term is identifiable, but that the full
measure $\nu$ is not, our aim is to examine where the boundary between what
can versus what cannot be identified lies. Toward this aim, one
direction to go is to
view~(\ref{I1}) as giving the first term of the expansion of the
``tail'' $F_{t}(\omega,[-u,u]^{c})$ near $0$,
and go further by assuming a series expansion such as
%
\begin{equation} \label{I2}%
F_{t}(\omega,[-u,u]^{c})\sim\sum_{i\geq1}\frac{a^{i}_{t}(\omega
)}{u^{\beta
_{i}}}\qquad \mbox{as }u\downarrow0
\end{equation}
(the precise assumption is given in Section~\ref{secSBG}), with successive
powers $\beta_{1}=\beta>\beta_{2}>\beta_{3}>\cdots.$ Those $\beta
_{i}$'s will
be the ``successive BG indices.'' This series
expansion can, for example, result from the superposition of processes with
different BG indices, in a model consisting of a sum of such processes.

The question then becomes one of identifying the successive terms in that
expansion. The main theoretical result of the paper, which is somehow
surprising, is as follows: the first index $\beta_{1}$ is always identifiable,
as we already knew, but the subsequent indices $\beta_{i}$ which are bigger
than $\beta_{1}/2$ are identifiable, whereas those smaller are not. An
intuition for this particular value of the ``identifiability
boundary'' is as follows: in view of~(\ref{I2}) the
estimation of the $\beta_{i}$'s can only be based on preliminary estimations
of $F_{t}(\omega,[-u,u]^{c})$, or of an integrated (in time) version of this,
for a sequence $u_{n}\rightarrow0$. It turns out that, even in idealized
circumstances, an estimation of $F_{t}(\omega,[-u_{n},u_{n}]^{c})$ or
of its
integrated version has a rate of convergence $u_{n}^{-\beta_{1}/2}$
(there is
a central limit theorem for this), so that any term contributing to
$F_{t}(\omega,[-u_{n},u_{n}]^{c})$ by an amount less than $u_{n}^{-\beta
_{1}/2}$ is fundamentally unreachable: we can only hope to estimate a further
coefficient $\beta_{i}$ if it leads to a number of increments greater than
$u_{n}$ (which is of order $u_{n}^{-\beta_{i}}$) that is larger than the
sampling error in the number of terms generated by the first coefficient,
implying that any $\beta_{i}<\beta_{1}/2$ cannot be identified. This shows
that there are limits to our ability to identify these successive
terms, even
in the unrealistic situation where the process is fully observed, and the
behavior of $\nu$ around $0$ is only partly identifiable.

When the identifiability conditions are satisfied, and when the process is
observed at discrete times with mesh $\Delta_{n}$, we will construct
estimators of the parameters which are consistent as $\Delta
_{n}\rightarrow0$,
and determine their rate of convergence, which we will see are slow. In the
case we have only two indices $\beta_{1}>\beta_{2}$ with $\beta
_{2}>\beta
_{1}/2$, we will further compare the rates of the estimators we
exhibit, which
are semiparametric, to the optimal rate achievable in a~corresponding
parametric sub-model (the sum of two stable processes, plus a drift and a
Brownian motion).%


%
\begin{figure}

\includegraphics{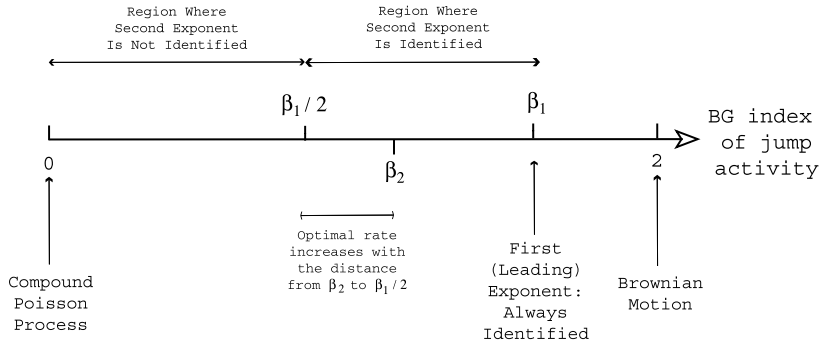}

\caption{Two BG component model: regions where the components are identified
versus not identified, and optimal rate of convergence.}%
\label{figtheorysummary}
\end{figure}

The main results of the paper are summarized in Figure~\ref{figtheorysummary}
for the two-component situation. We already noted that $\beta_{2}$ can be
identified only if it is bigger than $\beta_{1}/2;$ we will also see
that the
rate at which $\beta_{2}$ can be estimated increases as $\beta_{2}$ gets
closer to $\beta_{1}$, and conversely decreases as $\beta_{2}$ gets
closer to
$\beta_{1}/2$, in the limit dropping to $0$ as $\beta_{2}$ approaches
$\beta_{1}/2$, consistently with the loss of identification that occurs at
that point. Beyond the two-component model, we will provide general
identifiability conditions and rates of convergence for the leading and higher
order BG indices.

The paper is organized as follows. We first define the successive BG indices
in Section~\ref{secSBG}. In Section~\ref{secLC}, we study the
identifiability of the parameters appearing in the expansion, from a
theoretical viewpoint and in the special case of L\'evy processes. Then we
introduce consistent estimators for those parameters which we have
found to be
identifiable in the L\'evy case, hence proving de facto their
identifiability. This is done according to a two-step procedure, with
preliminary estimators given in Section~\ref{secD}, and final
estimators with
much faster rates in Section~\ref{secD2}. Unfortunately, although rates are
given, we were not able to show a central limit theorem for these estimators,
although such theorems ought to be available and would be crucial for
obtaining confidence bounds.

In principle, those estimators could be used on real data, but the
rates of
convergence for the higher order indices are, by necessity, quite slow. We
show in Section~\ref{secFisher} that the slow nature of these rates of
convergence is an inherent feature of the problem that cannot be improved
upon. This is perhaps not too surprising since the range of values of the
higher order indices that are identified is limited, and hence one would
expect the rate of convergence to deteriorate all the way to zero as one
approaches the region where identification disappears. We provide in Section
\ref{secMC} a simulation study for a~model featuring a stochastic volatility
plus two stable processes with different indices, the aim being to identify
these two indices, especially the higher order one. A realistic application
to high-frequency financial data, is out of the question for the typical
sample sizes that are currently available, but may be useful in the
future or
in different fields of applications where semimartingales are used and where
data are available in vast quantities, such as the study of Internet traffic
or turbulence data in meteorology. The results do also present theoretical
interest, especially as they set up bounds on what is asymptotically
identifiable in the jump measure of a semimartingale, and consequently
what is not.\vspace*{-1pt}

\section{The successive Blumenthal--Getoor indices}\label{secSBG}

Throughout the paper, $X$ is an It\^{o} semimartingale with characteristics
given by~(\ref{3}), on a filtered probability space $(\Omega,\mathcal{F}
,(\mathcal{F} _{t})_{t\geq0},\mathbb{P})$. The time horizon for the
observations is $T>0$, so the behavior of $X$ after time $T$ does not matter
for us below.

Our first aim is to give a precise meaning to an hypothesis like~(\ref{I2}).
Instead of requiring an expansion like this for all times $t$, we
rather use
the ``integrated version'' which uses the following family of (adapted,
continuous and increasing) processes:
%
\begin{equation}\label{30}%
u>0\quad\Rightarrow\quad\overline{A}(u)_{t} =\int_{0}^{t}F_{s}([-u,u]^{c})\,ds.
\end{equation}
The basic assumption is as follows:
%
\begin{assumption}
\label{AA1} There are a \textit{nonrandom} integer $j$, a strictly decreasing
sequence $(\beta_{i})_{1\leq i\leq j+1}$ of numbers in $[0,2)$ and a sequence
$(A^{i})_{1\leq i\leq j+1}$ of processes such that
%
\begin{equation}
\label{31}
t\in[0,T],\qquad 0<u\leq1\quad\Rightarrow\quad\Biggl\vert\overline{A}(u)_{t}
-\sum_{i=1}^{j}\frac{A_{t}^{i}}{u^{\beta_{i}}}\Biggr\vert\leq\frac
{A_{t}^{j+1}} {u^{\beta_{j+1}}}.
\end{equation}
Moreover, we have $A^{i}_{T}>0$ for $i=1,\ldots,j$.
\end{assumption}

If this assumption is satisfied with some $j\geq2$, it is also
satisfied with
any smaller integer. The processes $A^{i}$ and $A^{\prime}$ are nondecreasing
nonnegative, and they can always be chosen to be predictable.

Clearly, $\beta=\beta_{1}$ is the BG index, as introduced before, and the
following definition comes naturally in:
%
\begin{definition}
\label{DD1} Under Assumption~\ref{AA1}, the numbers $\beta_{1},\beta
_{2},\ldots, \beta_{j}$ are called the successive BG indices of the process
$X$ over the time interval $[0,T]$, and the variables $A^{i}_{T}$ are called
the associated integrated intensities.
\end{definition}
%
\begin{example}
\label{E21} Let $Y^{1},\ldots,Y^{j}$ be independent stable processes with
indices $\beta_{1}>\cdots>\beta_{j}$. Then $X=Y^{1}+\cdots+Y^{j}$ satisfies
(\ref{31}) with $A^{j+1}=0$\vadjust{\goodbreak} and the successive indices and integrated
intensities are $\beta_{i}$ and $Ta_{i}$, where $a_{i}=\lim_{u\to0}%
u^{\beta_{i}} F^{i}([-u,u]^{c})$, and $F^{i}$ is the L\'evy measure of
$Y^{i}$.

If the $Y^{i}$'s are tempered stable processes [see \citet{rosinski07}] the
same is true, provided $\beta_{j}>\beta_{1}-1$.\vspace*{-2pt}
\end{example}
%
\begin{example}
\label{E22} A semimartingale consisting of a continuous component and a jump
part driven by a sum of such processes also satisfies~(\ref{31}). Let
$X_{t}=X_{0}+Z_{t}+\sum_{i=1}^{j}\int_{0}^{t}H_{s}^{i}\,dY_{s}^{i}$, with
$Z$ a
continuous It\^{o} semimartingale and~$Y^{i}$ as in the previous
example and
$H^{i}$ locally bounded predictable processes with $\int_{0}^{T}|H^{i}%
_{s}|^{\beta_{i}}\,ds>0$. The successive BG indices are again the $\beta
_{i}%
$'s, with the associated integrated intensities
\[
A_{T}^{i}=a_{i}\int_{0}^{T}|H_{s}^{i}|^{\beta_{i}}\,ds.\vspace*{-2pt}
\]
\end{example}

\begin{remark}
\label{R19} We have taken a finite family of possible indices $\beta_{i}$.
Nothing prevents us from taking an infinite sequence: we simply have to assume
that Assumption~\ref{AA1} holds for all $j$, with additionally $\lim
_{i\to\infty} \beta_{i}=0$. However, in view of the restriction imposed on
the BG indices by our main theorems below about identifiability, this more
general situation has no statistical interest.\vspace*{-2pt}
\end{remark}
%
\begin{remark}
Assumption~\ref{DD1} imposes a certain structure on the behavior of the jump
measure of the process near $0$. It is important to note that it does not
restrict in any way the behavior of the jump measure away from~$0$. Although
most models used in practice and with infinite activity jumps satisfy this
assumption, the Gamma process does not: although it (barely) exhibits infinite
activity, its BG index is $0$, and $\overline{A}(u)_{t}$ is of order
$\log(1/u)$.\vspace*{-2pt}
\end{remark}

In Assumption~\ref{AA1}, expansion~(\ref{31}) is central, but one may
wonder about the additional requirement $A^{i}_{T}>0$. So, we end this section
with some comments and extensions, which may look complicated and are not
necessary for the rest of the paper, but which we think are useful and
somewhat enlightening.\vspace*{-2pt}

\begin{Extension}\label{exten1} In Assumption~\ref{AA1} positive and negative
jumps are treated in the same way. In practice, it might be useful for
modeling purposes to establish the behavior of positive and negative jumps
separately. Toward this end, one can replace~(\ref{30}) by
\[
\overline{A}(u)_{t}^{(+)} =\int_{0}^{t}F_{s}((u,\infty))\,ds,\qquad
\overline{A}(u)_{t}^{(+)} =\int_{0}^{t}F_{s}((-\infty,-u))\,ds.
\]
Then, if one is interested in positive jumps only, say, one replaces
(\ref{31}) by a~similar expansion for $\overline{A}(u)_{t}^{(+)}$: all the
content of the paper still holds, mutatis mutandis, under this modified
assumption, for positive jumps. The same is true of negative jumps, and the
``positive'' and ``negative'' successive BG indices can of course be
different.\vadjust{\goodbreak}
\end{Extension}
%
\begin{Extension}\label{exten2}
Now we come to the requirement $A^{i}_{T}>0$,
which in Assumption~\ref{AA1} is supposed to hold for all (or, almost all)
$\omega$. This is of course unlikely to hold for the terminal time $T$, unless
it holds for all $t>0$, and even unless the processes $A^{i}$ are strictly
increasing. In Example~\ref{E22}, this amounts to suppose that none of
processes $H^{i}$ vanishes. However, it might be relevant in practice
to allow
for each $H^{i}$ to vanish on some (possibly random) time intervals: we then
can have different components of the model turned on and off at
different times.

Thus, let us examine what happens if we relax the requirements $A^{i}_{T}>0$.
For any particular outcome $\omega$, the (first) BG index of the
process $X$
is $\beta_{i}$, where $i$ is the smallest integer such that
$A^{i}_{T}>0$, and
if all of them vanish one only knows that the BG index is not bigger than
$\beta_{j+1}$. The same applies to further indices. In other words, one can
define a partition of $\Omega$ indexed by all subsets $D$ of $\{1,\ldots
,j\}$
as follows:
%
\begin{equation} \label{37}%
\Omega_{T}(D)=\biggl( \bigcap_{i\in D}\{A_{T}^{i}>0\}\biggr) \cap\biggl(
\bigcap_{i\in\{1,\ldots,j\}\setminus D}\{A_{T}^{i}=0\}\biggr).
\end{equation}
Then, for any $\omega$, the successive BG indices of $X$ over $[0,T]$
and the
associated intensities are the numbers $\beta^{\prime}_{1}(\omega
),\ldots,\beta^{\prime}_{J}(\omega)$ and $\Gamma_{i}(\omega)$, defined as
%
\begin{eqnarray} \label{36}%
J(\omega)&=&m,\qquad \beta^{\prime}_{i}(\omega)=\beta_{l_{i}},\nonumber\\[-8pt]\\[-8pt]
\Gamma_{i}(\omega)&=&A_{t}^{l_{i}}(\omega) \qquad\mbox{if } \omega\in\Omega
_{t}(\{l_{1},\ldots,l_{m}\}).\nonumber
\end{eqnarray}
On the set $\Omega_{T}(\varnothing)$, which is not necessarily empty,
we have
$J=0$ and no~$\beta^{\prime}_{i}$'s.

All results of this paper are true if we relax $A^{i}_{T}>0$ in Assumption
\ref{AA1}, provided we replace $j$ by $J$ and the $\beta_{i}$'s by
the~$\beta^{\prime}_{i}$'s, \textit{in restriction to the set}~$\Omega_{T}(D)$:
this is indeed very easy, because on this set the process $X$ coincides
at all
times $t\in[0,T]$ with a process $X^{\prime}$ with satisfies Assumption
\ref{AA1} as stated above, with $(j,\beta_{1},\ldots,\beta_{j},\beta_{j+1})$
substituted with $(m,\beta_{l_{1}},\ldots,\beta_{j_{m}},\beta_{j+1})$, when
$D=\{l_{1},\ldots,l_{m}\}$.
\end{Extension}

\section{\texorpdfstring{Identifiability in the L\'{e}vy case}{Identifiability in the Levy case}}\label{secLC}

Loosely speaking, in an asymptotic statistical framework,
identifiability of a
parameter means the existence of a sequence of estimators which is (weakly)
consistent. Identifiability can be ``proved'' by exhibiting such a sequence.
It can be ``disproved'' by theoretical arguments, such as the fact that
if the
parameter is identifiable in our high-frequency observations setting, then,
were the path $t\mapsto X_{t}$ fully observed on $[0,T]$, it would enjoy
``nonasymptotic'' identifiability in the sense that its value is almost
surely known. For example, in the simple model $X_{t}=bt+W_{t}$ the parameter
$b$ does not enjoy this nonasymptotic property because the laws of the
process $X$ (restricted to $[0,T]$) are all equivalent when $b$ varies, and
thus $b$ is even less identifiable in the asymptotic setting.

Disproving identifiability is usually a hard task, especially in a
nonparametric setting. However, if a parameter is not identifiable for a
certain class of models, it is of course not identifiable for any wider class.

These arguments lead us to consider the very special situation of a L\'evy
processes $X$, with L\'{e}vy--Khintchine characteristics $(b,c,F)$ [see
(\ref{1})] when the path $t\mapsto X_{t}$ is fully observed on $[0,T]$. In
this section we are interested in nonasymptotic identifiability of those
characteristics, or functions of them. Note that, were $T$ infinite, the
triple $(b,c,F)$ would be identifiable because, for example, one would know
the values of all the i.i.d. increments $X_{n+1}-X_{n}$, giving us almost
surely the law of $X_{1}$, which in turn determines the triple $(b,c,F)$.

This is no longer the case when, as in this paper, the time interval $[0,T]$
is finite. In this case, we give a formal definition of
identifiability. We
use $Q_{b,c,F}$ to denote the law of the process $X$, restricted to the
interval $[0,T]$ ($T$ is kept fixed all throughout). So $Q_{b,c,F}$ is a
probability measure on the Skorokhod space $\mathbb{D}=\mathbb{D}
(|0,T],\mathbb{R})$. We also let $\mathcal{T}$ be some given subset of all
possible triples $(b,c,F)$.
%
\begin{definition}
\label{DD2} A function $H$ is \textit{identifiable on the class $\mathcal{T}$}
if, for any two $(b,c,F)$ and $(b^{\prime},c^{\prime},F^{\prime})$ in
$\mathcal{T}$ such that $H(b^{\prime},c^{\prime},F^{\prime})\neq
H(b,c,F)$, we
have $Q_{b,c,F}\perp Q_{b^{\prime},c^{\prime},F^{\prime}}$ (\textit{i.e., the
two measures $Q_{b,c,F}$ and $Q_{b^{\prime},c^{\prime},F^{\prime}}$ are
mutually singular}).
\end{definition}

The rationale behind this definition is as follows: if $H$ is
identifiable and
$(b,c,F)\in\mathcal{T}$, and $X$ is drawn according to the law $Q_{b,c,F}$,
then we can discard with probability $1$ any fixed $(b^{\prime
},c^{\prime
},F^{\prime})\in\mathcal{T}$ such that $H(b^{\prime},c^{\prime
},F^{\prime
})\neq H(b,c,F)$. Unfortunately, this does not mean that we can (almost
surely) reject all $(b^{\prime},c^{\prime},F^{\prime})$ with
$H(b^{\prime
},c^{\prime},F^{\prime}) \neq H(b,c,F)$ simultaneously: this stronger property
is (almost) never satisfied.

There exists a criterion for mutual singularity of $Q_{b,c,F}$ and
$Q_{b^{\prime},c^{\prime},F^{\prime}}$; see Remark IV.4.40 of
\citet{jacodshiryaev2003}. We have a Lebesgue decomposition $F^{\prime
}=f\bullet F+F^{\prime}{}^{\perp}$ of $F^{\prime}$ with respect to $F$, with
$f$ a nonnegative Borel function and $F^{\prime}{}^{\perp}$ a measure
supported by an $F$-null set. Then $Q_{b^{\prime},c^{\prime},F^{\prime
}}\perp
Q_{b,c,F}$ if and only if at least one of the following five properties is
violated:
%
\begin{equation}\label{LC1}
\cases{
\displaystyle F^{\prime}{}^{\perp}(\mathbb{R})<\infty,\vspace*{2pt}\cr
\displaystyle \alpha(F,F^{\prime})=\int\bigl( |f(x)-1|^{2}\wedge|f(x)-1|\bigr)
F(dx)<\infty,\vspace*{2pt}\cr
\displaystyle \alpha^{\prime}(F,F^{\prime})=\int_{\{|x|\leq1\}}|x| |f(x)-1|
F(dx)<\infty,
\vspace*{2pt}\cr
\displaystyle c=0\quad\Rightarrow\quad b^{\prime}=b-\int_{\{|x|\leq1\}}x \bigl(f(x)-1\bigr) F(dx),\cr
c^{\prime}=c.}
\end{equation}

It clearly follows that the function $H(b,c,F)=c$ is identifiable on
any class~$\mathcal{T}$ (a well-known fact). The function $H(b,c,F)=b$ is not
identifiable in general; however, on the class of all $(b,c,F)$ having $c=0$
and $\int_{\{|x|\leq1\}}|x|F(dx)<\infty$ the function $H(b,c,F)=\widehat
{b}=b-\int_{\{|x|\leq1\}}xF(dx)$ (which is the ``real'' drift, in the sense
that $X_{t}=\widehat{b}t+\sum_{s\leq t} \Delta X_{s}$) is identifiable.

In the sequel we are not interested in $b$ or $c$, but in $F$ only.
That is,
we are looking at functions $H=H(F)$. This leads us to consider classes
of the
form\looseness=-1
%
\begin{equation}\label{LC2}%
\mathcal{T}=\mathbb{R}\times\mathbb{R}_{+}\times\mathcal{T}_{3}
\qquad \mbox{where }\mathcal{T}_{3}\mbox{ is a set of L\'{e}vy measures.}
\end{equation}\looseness=0
In words, we want no restriction on the parameters $b$ and $c$. Of course
$\mathcal{T}_{3}$ should not be a singleton, and $H(F)$ should not be constant
on $\mathcal{T}_{3}$, otherwise the identifiability problem is empty.

The following example is clear:
%
\begin{example}
\label{ELC1} If $\mathcal{T}_{3}$ is a set of measures which coincide with
some given~$F$ on a neighborhood of $0$, then by~(\ref{LC1}) no nontrivial
$H(F)$ is identifiable on~$\mathcal{T}$.
\end{example}

This implies that, in the best-case scenario, a function $H(F)$ can be
identifiable only if it depends on the ``behavior of the
measure $F$ around $0$.'' Giving a necessary and
sufficient condition for identifiability of such a function, other than saying
that one of the properties in~(\ref{LC1}) fails when $H(F)\neq
H(F^{\prime})$,
seems out of reach. However, this is possible for some specific, but relatively
large, classes of sets $\mathcal{T}_{3}$, with a priori relatively
surprising results. Below we introduce such a class, in order to illustrate
the nature of the available results.
%
\begin{definition}[(The class $\mathcal{T}_{3}^{(1)}$ of L\'{e}vy measures)]
We say that a L\'{e}vy measure $F$ belongs to this class if we have
%
\begin{equation}\label{400}\quad
\cases{
\displaystyle F(dx) = \widetilde{F}(dx)+\sum_{i=1}^{\infty}\frac{a_{i}\beta_{i}
}{|x|^{1+\beta_{i}}} 1_{[-\eta,\eta]}(x)\,dx, \qquad \mbox{where $\eta>0$
and}\vspace*{2pt}\cr
\displaystyle\qquad \mbox{\hphantom{ii}(i)\quad}
0\leq\beta_{i+1}\leq\beta_{i}<2,\qquad\beta_{i}>0\quad\Rightarrow\quad
\beta_{i}>\beta_{i+1},\vspace*{2pt}\cr
\displaystyle\qquad \mbox{\hphantom{(iii)}\quad}\lim_{i\rightarrow\infty}\beta_{i}=0,\vspace*{2pt}\cr
\displaystyle\qquad \mbox{\hphantom{i}(ii)\quad} a_{i}>0\quad\Leftrightarrow\quad\beta_{i}>0,\vspace*{2pt}\cr
\displaystyle\qquad \mbox{(iii)\quad} 0<\sum_{i=1}^{\infty}a_{i}<\infty,\vspace*{2pt}\cr
\displaystyle\qquad \mbox{\hspace*{1pt}(iv)\quad} \widetilde{F} \mbox{ is a finite measure supported by
}[-\eta,,\eta]^{c}.}
\end{equation}
\end{definition}

Parts (i) and (ii) together ensure the uniqueness of the numbers $(a_{i},\beta_{i})$
in the representation of $F$, whereas if this representation holds for some
$\eta>0$, it also holds\vadjust{\goodbreak} for all $\eta^{\prime} \in(0,\eta)$, with the same
$(a_{i},\beta_{i})$. Part (iii) ensures that the infinite sum in the representation
converges, without being zero (so equivalently, $a_{1}>0$, or $\beta_{1}>0$).

The class $\mathcal{T}_{3}^{(1)}$ contains all sums of symmetric stable
L\'evy
measures. On the other hand, it is contained in the class of all L\'evy
measures $F$ of a L\'evy process satisfying Assumption~\ref{AA1}: the latter
is the class $\mathcal{T}_{3}^{(2)}$ of all $F$ such that
%
\begin{equation}
\label{4050}
u\leq1 \quad\Rightarrow\quad\Biggl\vert F([-u,u]^{c}) -\sum_{i=1}^{j}%
\frac{a_{i}}{u^{\beta_{i}}}\Biggr\vert\leq\frac{a^{\prime}} {u^{\beta
_{j+1}}}
\end{equation}
for $2>\beta_{1}>\cdots>\beta_{j+1}\geq0$ and $a_{i}>0$ for $i=1,\ldots
,j$ and
$a^{\prime}\geq0$, and those conditions are implied by~(\ref{400}), for any
$j\leq\sup(i\dvtx\beta_{i}>0)$, with the same $\beta_{i}$ and~$a_{i}$.

Considering $a_{i}$ and $\beta_{i}$ as functions on $\mathcal{T}_{3}^{(1)}$,
the identifiability result goes as follows:
%
\begin{theorem}
\label{TLC1} In the previous setting, the following holds:

\begin{longlist}
\item
The functions $\beta_{1}$ and $a_{1}$ are identifiable on the set
$\mathcal{T}^{(1)}_{3}$.

\item For any given $i\geq2$, the functions $\beta_{i}$ and $a_{i}$ are
identifiable on the subset $\mathcal{T}^{(1)}_{3}(i)=\{F\in\mathcal{T}
^{(1)}_{3}\dvtx\beta_{i}(F)\geq\beta_{1}(F)/2\}$ of $\mathcal
{T}^{(1)}_{3}$, and
they are not on the complement $\mathcal{T}^{(1)}_{3}\setminus\mathcal
{T}%
^{(1)}_{3}(i)$.
\end{longlist}
\end{theorem}
%
\begin{remark}
\label{R48} As mentioned in the ``first extension'' described in the previous
section, a similar statement is true if we replace the first line of
(\ref{400}) by
\[
F(dx)=\widetilde{F}(dx)+\sum_{j=1}^{\infty}\biggl(
\frac{a_{i}^{(+)}\beta_{i}^{(+)} }{|x|^{1+\beta_{i}}} 1_{(0,\eta]}%
(x)+\frac{a_{i}^{(-)}\beta_{i}^{(-)} }{|x|^{1+\beta_{i}}} 1_{(-\eta
,0)}(x)\biggr)\,dx
\]
with both families $(\beta_{i}^{(\pm)},a_{i}^{(\pm)})$ satisfying (i)--(iii).
Then the theorem above holds for both these families, with exactly the
same proof.
\end{remark}
%
\begin{remark}
\label{R50} As said before, any L\'{e}vy process $X$ whose L\'{e}vy
measure~$F$ is in $\mathcal{T}_{3}^{(1)}$ satisfies Assumption~\ref{AA1}, but the
converse is far from being true, so, even for L\'{e}vy processes, the
identifiability question is not completely solved under Assumption~\ref{AA1}.
More precisely, as the estimation results will show below,~(\ref{4050})
implies the ``positive'' identifiability
results [(i) and the first part of (ii) of Theorem~\ref{TLC1}] for L\'{e}vy
processes, but \textit{not} the ``negative''
results [second part of (ii)].

For example, consider the class $\mathcal{T}_{3}^{(3)}$ of all measure
of the
form
\[
F(dx)=\frac{a_{1} \beta_{1}}{x^{1+\beta_{1}}}
1_{(0,1]}(x)\,dx+G(dx)\qquad
\mbox{with }G = a_{2}\sum_{n\geq1}\varepsilon_{1/n^{1/\beta_{2}}}(dx)
\]
and $0<\beta_{2}<\beta_{1}<2$ and $a_{1},a_{2}>0$. Any such $F$ satisfies
(\ref{405}), but not~(\ref{400}). On $\mathcal{T}_{3}^{(3)}$, all four
parameters $\beta_{1},\beta_{2},a_{1},a_{2}$ are identifiable without the
restriction $\beta_{2}\geq\beta_{1}/2$. This is of course due to the
fact that
the measure $G$ is singular, and any two measures $G$ and $G^{\prime}$
of the
same type with $(\beta_{2},a_{2})\neq(\beta_{2}^{\prime},a_{2}^{\prime
})$ have
a Lebesgue decomposition $G^{\prime}=g\bullet G+G^{\prime\perp}$ with
$G^{\prime\perp}(\mathbb{R})=\infty$ when $\beta_{2}\neq\beta
_{2}^{\prime}$
and $\alpha(G,G^{\prime})=\infty$ when $\beta_{2}=\beta_{2}^{\prime}$ and
$a_{2}\neq a_{2}^{\prime}$.

We emphasize again that this example is quite singular, and verify here the
fairly general principle that the less regular a statistical problem
is, the
easier it is to solve in the sense that more parameters can be
estimated, and
often with faster rates.
\end{remark}
%
\begin{remark}
\label{R47} The class $\mathcal{T}_{3}^{(2)}$ may be bigger than
$\mathcal{T}_{3}^{(1)}$, but it is \textit{very far} from containing all
possible L\'evy measures. Indeed, any decreasing right-continuous
function $f$
on $(0,\infty)$ with $f(x)\to0$ as $x\to\infty$ and $f(x)\leq
K/x^{\alpha}$
for $x\in(0,1]$, for some constants $K>0$ and $\alpha\in(0,2)$, is the
symmetrical tail $f(x)=F([-x,x]^{c})$ of a L\'evy measure, although of course
it does not need to be equivalent to $a/x^{\beta}$ as $x\to0$ for some
$\beta\in(0,2)$ and $a>0$: so~(\ref{31}) may fail even with $j=1$.
\end{remark}

\section{Discretely observed semimartingales: Preliminary estimators}\label{secD}

Now we turn to the more general case of semimartingales. The process
$X$ is
observed at the times $i\Delta_{n}$ for $i=0,1,\ldots,[T/\Delta_{n}]$ (where
$[x]$ denotes the integer part of the real~$x$). We thus observe the
increments
%
\begin{equation} \label{D1}%
\Delta_{i}^{n}X=X_{i\Delta_{n}}-X_{(i-1)\Delta_{n}}.
\end{equation}

The BG indices describes some properties of jumps, which are not observed.
However, when an increment $\Delta_{i}^{n}X$ is relatively large, say bigger
than $u_{n}$ with $u_{n}\gg\sqrt{\Delta_{n}}$, it is likely to be due to
jumps because the drift plus the continuous martingale part have
increments of
order of magnitude $\sqrt{\Delta_{n}}$. Moreover it turns out that it is
usually due to a single ``large'' jump of size bigger than $u_{n}$, although
of course the observed value $\Delta^{n}_{i}X$ is not exactly the jump size.
So one may expect the number of jumps with size bigger than $u_{n}=u$, over
the time interval $[0,t]$, to be the following number, or be relatively close
to it:
%
\begin{equation}\label{D2}%
U(u,\Delta_{n})_{t}=\sum_{i=1}^{[t/\Delta_{n}]}1_{\{\Delta_{i}^{n}X>u\}}.
\end{equation}
In order for the previous statement to actually be true, we need some
additional assumptions, though. Those are given in the following:
%
\begin{assumption}
\label{A2} The process $X$ is an It\^{o} semimartingale, and:

\begin{longlist}[(a)]
\item[(a)]
The processes $b_{t}$, $c_{t}$ are locally bounded.\vadjust{\goodbreak}

\item[(b)] We have Assumption~\ref{AA1} with $A_{t}^{i}=\int_{0}^{t}
a_{s}^{i}\,ds$ for
$i=1,\ldots,j+1$, where the processes $a^{i}$ are locally bounded.

\item[(c)] We have $\beta_{j}>\beta_{1}/2$.
\end{longlist}
\end{assumption}

Assumption~\ref{A2}(c) above may look strange, or too strong. However, in view
of the
identifiability results of the previous section, we cannot estimate
consistently $\beta_{i}$ if it is strictly smaller than $\beta_{1}/2$,
and as
a matter of fact, the estimators described below are consistent only if
$\beta_{i}>\beta_{1}/2$. Hence, since Assumption~\ref{AA1} for $j$
implies the
same for all $j^{\prime}<j$, (c) above is really \textit{not} a restriction,
but amounts to replacing $j$ in this assumption by $j\wedge\sup\{
i\dvtx\beta_{i}>\beta_{1}/2\} $.

Apart from (c), this assumption is satisfied in Examples~\ref{E21} and
\ref{E22}, and also by any L\'evy process satisfying Assumption~\ref{AA1}.

The estimation procedure is a two-step procedure, and in this section we
describe the first---preliminary---estimators. These estimators will be
consistent, but with very slow rates of convergence. This is why, in
the next
subsection, we will derive final estimators which exhibit much faster
(although still slow) rates.

Those preliminary estimators require the knowledge of a number
$\varepsilon>0$
which satisfies
%
\begin{equation}\label{118}%
i=1,\ldots,j-1\quad\Rightarrow\quad\beta_{i}-\beta_{i+1}\geq\varepsilon.
\end{equation}
Such an $\varepsilon$ always exists, but here we suppose that it is known,
somewhat in contradiction with the fact that the $\beta_{i}$ are
unknown. It
it is obviously quite difficult to estimate properly two contiguous indices
$\beta_{i}$ and $\beta_{i+1}$ when they are very close to to one
another. So
from a statistical viewpoint, the assumption $\beta_{i}-\beta_{i+1}%
>\varepsilon$ for some fixed $\varepsilon>0$ is natural. Moreover,
since we do
not know a priori which $\omega$ is observed, this amounts to
supposing that all possible values of the BG indices in the model
satisfy this
restriction. For models used in practice, this is not really a restriction
since these models rely on at most a small number of indices that are
separated from one another.

The key ingredient for constructing the estimators is the counting process
defined in~(\ref{D2}), evaluated at the terminal time $T$ and for suitable
values of~$u$. In particular, we choose a sequence $u_{n}$ satisfying
%
\begin{eqnarray}\label{D31}
u_{n}\rightarrow0, \qquad\Delta_{n}^{\rho}\leq Ku_{n}\hspace*{40pt}\nonumber\\[-8pt]\\[-8pt]
&&\eqntext{\mbox{with }%
\displaystyle \rho<\frac{1}{2+\beta_{1}}\wedge\frac{2}{\beta_{1}(3+\beta_{1})}%
\wedge\frac{4}{\beta_{1}(5+3\beta_{1})}.}
\end{eqnarray}
Of course $\rho>0$ above (otherwise $u_{n}\rightarrow0$ would fail). The
infimum of the upper bound for $\rho$ over all $\beta_{1}<2$ is $2/11$.
Therefore, since we do not a~priori know the values of
$\beta_{1}$, whereas as we will see the rates improve when the sequence
$u_{n}$ becomes smaller (termwise), it is thus advisable to take $\rho
=2/11$ above.\vadjust{\goodbreak}

The first-step estimation is done by induction on $i$. We choose $\gamma>1$,
and the estimators for $\beta_{1}$ and $A_{T}^{1}$ are
%
\begin{eqnarray}
\label{D3}\qquad
\widetilde{\beta}_{1}^{n}&=&\cases{\displaystyle
\frac{\log(U(u_{n},\Delta_{n})_{T}/U(\gamma u_{n},\Delta_{n})_{T})}{\log
\gamma}, &\quad if $U(\gamma u_{n},\Delta_{n})_{T}>0$,\vspace*{2pt}\cr
-1, &\quad otherwise,}
\nonumber\\[-8pt]\\[-8pt]
\widetilde{\Gamma}_{i}^{n}&=&(u_{n})^{\widetilde{\beta}_{1}^{n}}
U(u_{n},\Delta_{n})_{T}.%
\nonumber
\end{eqnarray}

For constructing the subsequent estimators, and with $\varepsilon$ in
(\ref{118}), we set
%
\begin{equation} \label{S10}%
u_{n,i}=u_{n}^{(\varepsilon/2)^{i-1}}
\end{equation}
(so $u_{n,1}=u_{n}$). We denote by $I(k,l)$ the set of all subsets of
$\{1,\ldots,k\}$ having~$l$ elements. Assuming that we know $\widehat
{\beta
}_{i}^{n}$ and $\widehat{\Gamma}_{i}^{n}$ for $i=1,\ldots,k-1$, for some
$k\in\{2,\ldots,j\}$, we set
%
\begin{eqnarray}\label{D4}
x&\geq&1\quad\Rightarrow\quad U^{n}(k,x)=\sum_{l=0}^{k-1}(-1)^{l} U(x \gamma
^{l} u_{n,k},\Delta_{n})_{T}\sum_{J\in I(k-1,l)}\gamma^{ \sum_{i\in
J}\widetilde{\beta}_{i}^{n}},\hspace*{-35pt}\nonumber\\
\widetilde{\beta}_{k}^{n}&=&\cases{\displaystyle
\frac{\log( U^{n}(k,1)/U^{n}(k,\gamma)) }{\log(\gamma)}, &\quad
if $U^{n}(k,1)>0,U^{n}(k,\gamma)>0$,\vspace*{2pt}\cr
-1, &\quad otherwise,}\hspace*{-35pt}
\\
\widetilde{\Gamma}_{k}^{n}&=&u_{n,k}^{\widetilde{\beta}_{k}^{n}} \Biggl(
U(u_{n,k},\Delta_{n})_{T}-\sum_{l=1}^{k-1}\widetilde{\Gamma}_{l}^{n}
u_{n,k}^{-\widetilde{\beta}_{l}^{n}}\Biggr) .\nonumber\hspace*{-35pt}
\end{eqnarray}

Finally, in order to state the result, we need a further notation, for
$i=1,\ldots,j-1$ \mbox{(so when $j=1$ the following is empty)}:
%
\begin{equation} \label{S12}%
H_{i}=\frac{A_{T}^{i+1}}{A_{T}^{i} \log\gamma}\frac{\prod_{l=1}
^{i}( \gamma^{\beta_{l}-\beta_{i+1}}-1) }{\prod_{l=1}
^{i-1}( \gamma^{\beta_{l}-\beta_{i}}-1) }.
\end{equation}

\begin{theorem}
\label{TD1} Under Assumption~\ref{A2} and~(\ref{118}), for all
$i=1,\ldots
,j-1$ such that $\beta_{i+1}>\beta_{1}/2$, we have
%
\begin{equation} \label{D5}
\frac{\widetilde{\beta}_{i}^{n}-\beta_{i}}{u_{n,i}^{\beta_{i}-\beta_{i+1}}
}\stackrel{\mathbb{P}}{\longrightarrow}-H_{i},\qquad \frac{\widetilde{\Gamma}
_{i}^{n}-A_{T}^{i}}{u_{n,i}^{\beta_{i}-\beta_{i+1}} \log(1/u_{n,i}
)}\stackrel{\mathbb{P}}{\longrightarrow}\Gamma_{i} H_{i}.
\end{equation}
Moreover if $\eta=\beta_{j}-\beta_{j+1}\vee\frac{\beta_{1}}2>0$, the following
variables are bounded in probability:
%
\begin{equation}
\label{D6}
\frac{\widetilde{\beta}_{j}^{n}-\beta_{i}}{u_{n,j}^{\eta}}
,\qquad \frac{\widetilde{\Gamma}_{j}^{n}-A_{T}^{j}}{u_{n,j}^{\eta}
\log(1/u_{n,j})}.
\end{equation}
\end{theorem}

The estimator $\widetilde{\beta}_{1}^{n}$ is exactly the estimator
proposed in
\citet{yacjacod09b} for the leading BG index $\beta_{1}$. So, not only
does it
satisfy~(\ref{D5}) when $j\geq2$ or the tightness of~(\ref{D6}) when $j=1$,
but it also enjoys a central limit theorem centered at $\beta_{1}$ and with
rate $u_{n}^{\beta_{1}/2}$ as soon as $\beta_{2}<\beta_{1}/2$ (this property
implies $j=1$ here). Moreover, in this case one could prove
that~$\widetilde{\Gamma}_{1}^{n}$ also satisfies a CLT with the rate $u_{n}%
^{\beta_{1}/2} \log(1/u_{n})$, although we will not prove it, since the
emphasis here is on the case of several BG indices.\looseness=1

Some remarks are in order here:
%
\begin{remark}
\label{RS0} It is possible for the estimator $\widetilde{\Gamma
}^{n}_{i}$ to
be negative, in which case we may replace it by $0$, or by any other positive
number. It may also happen that the sequence $\widetilde{\beta
}^{n}_{i}$ is
not decreasing, and we can then reorder the whole family as to obtain a
decreasing family (we relabel the estimators of $A_{T}^{i}$
accordingly, of
course). All these modifications are asymptotically immaterial.
\end{remark}
%
\begin{remark}
\label{RS1} As mentioned in the Extension~\ref{exten2} at the end of Section
\ref{secSBG}, we can relax $A^{i}_{T}>0$ in Assumption~\ref{AA1}. Then the
above theorem is still valid, in restriction to the set $\Omega_{T}%
(\{l_{1},\ldots,l_{m}\})$ of~(\ref{37}), as soon as $\beta_{l_{m}}%
>\beta_{l_{1}}/2$.
\end{remark}
%
\begin{remark}
\label{RS2} Suppose that $j\geq2$. The limits in~(\ref{D5}) are pure bias,
hence precluding the existence of a proper\vspace*{1pt} central limit theorem. Note that
$H_{i}>0$ if $i<j$, so the bias for $\widehat{\beta}_{i}^{n}$ and for
$\widehat{\Gamma}_{i}^{n}$ are always negative and positive, respectively.

Note also that the rate of convergence for estimating $\beta_{i}$ when
$i\leq j-1$, say, is $u^{\beta_{i}-\beta_{i+1}}_{n,i}$, that is
$u_{n}^{(\beta _{i}-\beta_{i+1})(\varepsilon/2)^{i-1}}$. This is
exceedingly small, indeed. For example, suppose that we have three
indices $\beta_{1}>\beta _{2}>\beta _{3}>\frac{\beta_{1}}2$.
Then~(\ref{118}) implies necessarily $\varepsilon <\frac{\beta_{1}}2$, so
the best\vspace*{1pt} possible rate for $i=2$ would be less than, but
close to, $u_{n}^{(\beta_{2}-\beta_{3})\beta_{1}/4}$, upon taking
$\varepsilon$ close to $\frac{\beta_{1}}2$, which is of course
impossible because we do not know $\beta_{1}$ to start with.

In the previous example, if we suspect that $\beta_{1}$ is bigger than $1$,
say, it becomes (perhaps) not totally unreasonable to choose
$\varepsilon
=0.1$; the rates for $i=2$ and $i=3$ thus become $u_{n}^{(\beta
_{2}-\beta
_{3})/10}$ and $u_{n}^{(\beta_{3}-\beta_{1}/2)/100}$. This is of course
on top
of the fact that, because of~(\ref{D31}), $u_{n}$ is of order of magnitude
$\Delta_{n}^{2/11}$, by a conservative choice of $\rho$.
\end{remark}

\textit{Practical considerations}. Letting aside the slow convergence
rates, the previous result suffers from two main drawbacks:

(1) It requires to know the number of indices to be estimated (this is
implicit in Assumption~\ref{A2}).

(2) It requires to know a number $\varepsilon>0$
satisfying~(\ref{118}).\vadjust{\goodbreak}

About the first problem above, in real world one does not know the
number of indices. On the other hand, if Assumption~\ref{AA1} holds,
it seems
reasonable to suppose that it holds for all $j$, whereas the estimation is
made for those $\beta_{i}$ which are bigger than $\beta_{1}/2$ only. In
connection with this, we assume $\beta_{i}-\beta_{i+1}\geq\varepsilon$
for all
$i\leq j:= \sup(k\dvtx \beta_{k}>\beta_{1}/2)$, plus the property $\beta_{j}%
>\beta_{1}/2+ \varepsilon$. Then the aim becomes to estimate $\beta
_{i}$ and
$A^{i}_{T}$ for all $i\leq j$, with $j$ unknown.

Since the estimation procedure is done by induction on the successive indices,
one can start the induction as described above, and stop it at the
first $i$
such that $\widetilde{\beta}_{i}\leq\varepsilon+\widetilde{\beta}_{1}/2$.
Asymptotically, this procedure will deliver the ``correct'' answer (the proof
of this fact, not given below, is a simple extension of the proof of the
second claim of the theorem). In practice, however, the solution to this
stopping problem is not quite clear, since in particular the estimated
sequence $\widetilde{\beta}_{i}$ is not necessarily decreasing,
although it is
so asymptotically.

Problem 2 above is clearly more annoying. We have to admit that, in the
setting presented here, we have no theoretical solution for solving it. A
possible way out would be to make the estimation with several values of
$\varepsilon$, going downward, until the estimated differences
$\widetilde
{\beta}_{i}-\widetilde{\beta}_{i-1}$ all become significantly bigger
than the
chosen~$\varepsilon$, but no mathematical result so far is available in this
direction. In addition, since rates are very slow, the probability that
such a
difference is bigger than $\varepsilon$ when the true values satisfy
the same
inequality may be not close to $1$ (for finite, but even large, samples).

Nonetheless, bad as it looks, this condition is probably relatively innocuous
in practice: indeed, when two successive indices are very close to each
other, they are obviously very difficult to tell apart. So the problem is
practically meaningful only if the indices are a small number (as $2$,
$3$ or
perhaps~$4$) and reasonably well separated. Hence taking $\varepsilon
=0.1$ for
instance, as in Remark~\ref{R47}, seems to be safe enough.

\section{Discretely observed semimartingales: An improved method}\label{secD2}

The observation scheme is the same as in the previous section: $X$ is observed
at the times $i\Delta_{n}$ smaller or equal to some fixed terminal time $T$.

As already mentioned, the previous estimators converge at a \textit{very} slow
rate, especially for higher order indices; see Remark~\ref{RS2}. So, in order
to implement the estimation with any kind of reasonable accuracy, it is
absolutely necessary to come up with better estimators.

This is the aim of this section. Assuming Assumption~\ref{A2}, we also suppose
that we can construct preliminary estimators, such as in the previous section.
Exactly as there, we must know the number $j$ of BG indices that are to
be estimated.

The method consists in minimizing, at each stage $n$, a suitably chosen
contrast function $\Phi_{n}$. First we take an integer $L\geq2j$ and numbers
$1=v_{1}<v_{2}<\cdots<v_{L}$. We also choose positive weights $w_{k}$\vadjust{\goodbreak}
(typically $w_{k}=1$, but any choice is indeed possible), and we pick
truncation levels $u_{n}$ satisfying~(\ref{D31}). We also let $D$ be
the set
of all $(x_{i},\gamma_{i})_{1\leq i\leq j}$ with $0\leq x_{j}\leq x_{j-1}
\leq\cdots\leq x_{1}\leq2$ and $\gamma_{i}\geq0$. Then the contrast function
is defined on $D$ by
%
\begin{equation} \label{D2-2}%
\Phi_{n}(x_{1},\gamma_{1},\ldots,x_{j},\gamma_{j})=\sum_{l=1}^{L}
w_{l}\Biggl( U(v_{l}u_{n},\Delta_{n})_{T}-\sum_{i=1}^{j}\frac{\gamma_{i}
}{(v_{l}u_{n})^{x_{i}}}\Biggr) ^{2},
\end{equation}
where the sequence $u_{n}$ satisfies~(\ref{D31}). Then the estimation
goes as follows:

\begin{Step}\label{step1}
We construct preliminary estimators $\widetilde{\beta}^{n}_{i}$
(decreasing in $i$) and~$\widetilde{\Gamma}^{n}_{i}$ (nonnegative) for
$\beta_{i}$ and $A_{T}^{i}$ for $i=1,\ldots,j$, such that $(\widetilde
{\beta
}_{i}-\beta_{i}) /u_{n}^{\eta}$ and $(\widetilde{\Gamma}_{i}-A_{T}^{i}
)/u_{n}^{\eta}$ go to $0$ in probability for some $\eta>0$. For
example, we may
choose those described in the previous section (see Remark~\ref{RS0}): the
consistency requirement is fulfilled for any $\eta<(\varepsilon/2)^{j}$.
\end{Step}
\begin{Step}\label{step2}
We denote\vspace*{1pt} by $D_{n}$ the (compact and nonempty) random subset
of $D$
defined by $D_{n}=\{(x_{i},\gamma_{i})\in D\dvtx |x_{i}-\widetilde{\beta}^{n}
_{i}|\leq\alpha u_{n}^{\eta},|\gamma_{i}-\widetilde{\Gamma}^{n}_{i}
|\leq\alpha u_{n}^{\eta},\forall i=1,\ldots,j\}$, for some arbitrary (fixed)
$\alpha>0$. Then the final estimators $\overline{\beta}{}^{n}_{i}$
and~$\overline{\Gamma}{}^{n}_{i}$ will be
%
\begin{equation}
\label{D2-3}
(\overline{\beta}{}^{n}_{i},\overline{\Gamma}{}^{n}_{i})_{1\leq
i\leq
j}=\mathop{\arg\min}_{D_{n}}\Phi_{n}(x_{1},\gamma_{1},\ldots
,x_{j},\gamma_{j}).
\end{equation}
\end{Step}

\begin{theorem}
\label{TD2} Under Assumption~\ref{A2}, and for all choice of $v_{2}
,\ldots,v_{L}$ outside a $\lambda_{L-1}$-null set (depending on the
$\beta
_{i}$'s; $\lambda_{l}$ is the $l$-dimensional Lebesgue measure), the
sequences
%
\begin{equation}
\label{D2-4}
\frac{\overline{\beta}{}^{n}_{i}-\beta_{i}}{u_{n}^{\beta_{i}%
-\beta_{1}/2-\mu}},\qquad \frac{\overline{\Gamma}{}^{n}_{i}-\Gamma_{i}}%
{u_{n}^{\beta_{i}-\beta_{1}/2-\mu}}%
\end{equation}
are bounded in probability for all $i=1,\ldots,j$ and all $\mu>0$.
\end{theorem}

The rates obtained here are much faster than in Theorem~\ref{TD1}: we replace
$u_{n,i}^{\beta_{i}-\beta_{i+1}\vee(\beta_{1}/2)}$ by $u_{n}^{\beta_{i}%
-\beta_{1}/2}$, for two reasons: the exponent $\beta_{i}-\beta_{i+1}$ is
bigger than $\beta_{i}-\beta_{i+1}\vee(\beta_{1}/2)$, unless $i=j$; more
importantly, we replace the auxiliary truncation levels $u_{n,i}$ of
(\ref{S10})
by the original sequence $u_{n}$, which is much smaller when $i\geq2$, and
only subject to~(\ref{D31}). We will examine in the next section how
far from
optimality those rates are.
%
\begin{remark}
\label{R40} As stated, and as seen from the proof, we only need $L=2j$,
and choosing $L>2j$ does not improve the asymptotic properties.
However, from a practical viewpoint, it is probably wise to take $L$
bigger than $2j$ in order to smooth out the contrast function somehow,
especially for (relatively) small samples. A~choice of the weights
$w_{l}>0$ other than $w_{l}=1$, such as~$w_{l}$ decreasing in $l$,
may
serve to put less emphasis on the large truncation values $u_{n}v_{l}$
for which less data are effectively used.\vadjust{\goodbreak}
\end{remark}
%
\begin{remark}
\label{R41} The result does not hold (or at least we could not prove
it) for
\textit{all} choices of the $v_{l}$'s, but only when $(v_{2},\ldots,v_{L})$
(recall $v_{1}=1$) does not belong to some Lebesgue-null set $G(\beta
_{1},\ldots,\beta_{j})$. This seems a priori a~serious restriction, because
$(\beta_{1},\ldots,\beta_{j})$ is unknown. In practice, we choose a priori
$(v_{2},\ldots,v_{L})$, so we may have bad luck, just as we may have
bad luck
for the outcome $\omega$ which is drawn$\ldots.$

We may also do the estimation for a number of different choices for the
weights and/or values of $L\geq2j$ and compare or average the results. This
should contribute to weaken the numerical instability inherent to minimization
problems such as~(\ref{D2-3}). This numerical instability is similar to the
one occurring in nonlinear regression problems.

We have to state, however, that these problems, just as those stated in the
``practical considerations'' of the previous section, are not fully addressed
in this paper, and they are probably quite difficult to overcome. Our emphasis
here is more on theoretical results, and on the possibility of
performing the
estimation with reasonable rates (see, however, Section~\ref{secMC}
below, to
see how the problem of finding a ``good'' $\varepsilon$ and doing preliminary
estimation in our simulation study is skipped, without affecting the quality
of the procedure in any noticeable way).\vspace*{-1pt}
\end{remark}

\section{Optimality in a special case}\label{secFisher}\vspace*{-1pt}

\subsection{Why the convergence rates are necessarily slow}

Intuitively, the fact that we are right at the boundary between
identifiability and lack thereof suggests that we should expect the
rate, as
we approach the loss of identifiability boundary, to deteriorate all
the way
to zero. In order to quantify precisely how slow the rates of
convergence for
the estimators of the second (and higher) index must be, even in ideal
circumstances, we study a simple parametric model of the following
form. Let
$W$ be a Brownian motion and $Y^{1},Y^{2}$ be two independent standard
symmetric stable processes, and set
%
\begin{equation} \label{eqFisher+2Stable}%
X_{t}=bt+\sigma W_{t}+Y_{t}^{1}+Y_{t}^{2}.
\end{equation}

Each $Y^{i}$ depends on two parameters, the index $\beta_{i}$ and a scale
parameter~$a_{i}$, the latter being characterized by the fact that the
L\'evy
measure of~$Y^{i}$ is
%
\begin{equation}
\label{LR4-7}
F^{j}(dx)=\frac{a_{j} \beta_{j}}{|x|^{1+\beta_{j}}}\,dx.
\end{equation}
We have six parameters,
%
\begin{equation}
\label{IJ2}
b\in\mathbb{R},\qquad c=\sigma^{2}>0,\qquad a_{1},a_{2}%
>0,\qquad 0<\beta_{2}<\beta_{1}<2,
\end{equation}
among which $b$ is not identifiable, and $c,\beta_{1},a_{1}$ are identifiable,
and $(\beta_{2},a_{2})$ are identifiable if and only if $\beta_{2}\geq
\beta_{1}/2$. In what follows, we restrict our attention to the four
parameters $\beta_{1},\beta_{2},a_{1},a_{2}$.

In order to find at which rate it is possible to estimate these four
parameters, when $X$ is observed\vadjust{\goodbreak} at the discrete times $(i\Delta
_{n}\dvtx i=0,1,\ldots,[T/\Delta_{n}])$ and $\Delta_{n} \rightarrow0$, we
study the
behavior of the Fisher information matrix. Due to the fact that $X$ is a
L\'{e}vy process, the information matrix at stage $n$ is $[T/\Delta_{n}]$
times the information matrix obtained when we observe only the
variable~$X_{\Delta_{n}}$; since the variable $X_{\Delta}$ admits a density
$x\mapsto
p(_{\Delta}(x|c,\beta_{1},a_{1},\allowbreak \beta_{2},a_{2})$ which is $C^{\infty}$ in
$x$, and also in $(c,\beta_{1},a_{1},\beta_{2},a_{2})$ on the domain defined
by~(\ref{IJ2}), it is no wonder that Fisher's information $I_{\Delta}$
for a
single observation~$X_{\Delta}$ (recall $X_{0}=0$) exists, and we can study
its behavior as $\Delta\rightarrow0$.

Only the diagonal entries are important for the various rates of convergence,
so we only need to focus on the following diagonal entries of this matrix:
\[
I_{\Delta}^{\beta_{1}\beta_{1}},I_{\Delta}^{a_{1}a_{1}},I_{\Delta}%
^{\beta_{2}\beta_{2}}, I_{\Delta}^{a_{2}a_{2}}.
\]

The main result of this section follows, giving the asymptotic order of the
relevant terms in Fisher's information:
%
\begin{theorem}
\label{theo-fisher}We have the following equivalences, as $\Delta
\rightarrow
0$:
\begin{eqnarray*}
I_{\Delta}^{\beta_{1}\beta_{1}} & \sim & \frac{a_{1} }{2(2-\beta_{1}%
)^{\beta_{1}/2} c^{\beta_{1}/2}} \Delta^{1-\beta_{1}/2} (\log
(1/\Delta))^{2-\beta_{1}/2},\\
I_{\Delta}^{a_{1}a_{1}} & \sim & \frac{2\beta_{1}c_{\beta_{1}} a_{1}
^{\beta_{1}}}{(2-\beta_{1})^{\beta_{1}/2} \sigma^{\beta_{1}} a_{1}^{2}
} \frac{\Delta^{1-\beta_{1}/2}}{(\log(1/\Delta))^{\beta_{1}/2}}%
\end{eqnarray*}
and also, provided $\beta_{2}>\beta_{1}/2$,
\begin{eqnarray*}
I_{\Delta}^{\beta_{2}\beta_{2}} & \sim &\frac{a_{2}^{2} \beta_{2}^{2}}
{2a_{1} \beta_{1}(2 \beta_{2}-\beta_{1})(2-\beta_{1})^{\beta_{2}-\beta_{1}/2}
c^{\beta_{2}-\beta_{1}/2}}\\
&&{}\times \Delta^{1-\beta_{2}+\beta_{1}/2} (\log
(1/\Delta))^{2-\beta_{2}+\beta_{1}/2},\\
I_{\Delta}^{a_{2}a_{2}} & \sim & \frac{2\beta_{2}^{2}} {a_{1} \beta
_{1}(2\beta_{2}-\beta_{1})(2- \beta_{1})^{\beta_{2}-\beta_{1}/2}
c^{\beta
_{2}-\beta_{1}/2}} \frac{ \Delta^{1-\beta_{2}+\beta_{1}/2}}{(\log
(1/\Delta))^{\beta_{2}-\beta_{1}/2}}.
\end{eqnarray*}
\end{theorem}
%
\begin{remark}
We are not concerned here with the identification and estimation of the
volatility parameter $c$; the term $I_{\Delta}^{cc}$ in a simpler model
has been studied in \citet{yacjacod08}, as well as
$I_{\Delta}^{a_{1}a_{1}}$ when $a_{2}=0$ (i.e., when there is only one
stable process on top of the Brownian motion). The asymptotic
equivalent for the term $I_{\Delta}^{a_{1}a_{1}}$ of course reduces to
(4.11) of that paper, with $\alpha=\beta_{1}$, $\beta=2$,
$\theta=a_{1}$, up to a change of parametrization for $a_{1}$, since
here we use the parametrization~(\ref{LR4-7}) which corresponds to the
notation of Assumption~\ref{AA1}, which is fulfilled here.
\end{remark}

Coming back to the original problem, we deduce that it should be
possible in
principle to find estimators $\widehat{\beta}_{i}^{n}$ and $\widehat{a}
_{i}^{n}$ having the following properties:
%
\begin{eqnarray}
\label{119}
\frac{(\log(1/\Delta_{n}))^{1-\beta_{1}/4}}{\Delta_{n}^{\beta_{1}/4}
} (\widehat{\beta}_{1}^{n}-\beta_{1}) & \stackrel{\mathcal
{L}}{\longrightarrow
} & \mathcal{N}(0,1/T\mathcal{I}^{\beta_{1}\beta_{1} }),\nonumber\\
\frac{1}{\Delta_{n}^{\beta_{1}/4} (\log(1/\Delta_{n}))^{\beta_{1}/4}
} (\widehat{a}_{1}^{n}-a_{1}) & \stackrel{\mathcal{L}}{\longrightarrow} &
\mathcal{N}(0,1/T\mathcal{I}^{a_{1}a_{1}}),\nonumber\\[-8pt]\\[-8pt]
\frac{(\log(1/\Delta_{n}))^{1-\beta_{2}/2+\beta_{1}/4}}{\Delta
_{n}^{\beta
_{2}/2-\beta_{1}/4}} (\widehat{\beta}_{2}^{n}-\beta_{2}) & \stackrel
{\mathcal{L}}{\longrightarrow} & \mathcal{N}(0,1/T\mathcal{I}^{\beta_{2}
\beta_{2}}),\nonumber\\
\frac{1}{\Delta_{n}^{\beta_{2}/2-\beta_{1}/4} (\log(1/\Delta
_{n}))^{\beta
_{2}/2-\beta_{1}/4}} (\widehat{a}_{2}^{n}-a_{2}) & \stackrel{\mathcal{L}
}{\longrightarrow} & \mathcal{N}(0,1/T\mathcal{I}^{a_{2}a_{2}}),
\nonumber
\end{eqnarray}
where $\mathcal{I}^{\beta_{1}\beta_{1}}$, $\mathcal{I}^{a_{1}a_{1}}$,
$\mathcal{I}^{\beta_{2}\beta_{2}}$ and $\mathcal{I}^{a_{2}a_{2}}$ are the
constants in front of the term involving $\Delta$ in the equivalences above,
for $I_{\Delta}^{\beta_{1}\beta_{1}}$, $I_{\Delta}^{a_{1}a_{1}}$,
$I_{\Delta
}^{\beta_{2}\beta_{2}}$ and $I_{\Delta}^{a_{2}a_{2}}$, respectively.
Conversely, by the Cram\'er--Rao lower bound, Theorem~\ref{theo-fisher} also
implies that it will be impossible to find consistent estimators with faster
rates of convergence, or smaller asymptotic variance, that those
exhibited in
(\ref{119}).

Note that these rates are consistent with the results of Theorem~\ref{TLC1}.
The first two convergences above shows that it is always possible to estimate
consistently $\beta_{1}$ and $a_{1}$, the third one implies consistency for
$\beta_{2}$ only if $\beta_{2}\geq\beta_{1}/2$, and the last one implies
consistency for $a_{2}$ only if $\beta_{2}>\beta_{1}/2$. The last statement
seems contradictory with Theorem~\ref{TLC1} when $\beta_{2}=\beta
_{1}/2$, but
of course it is possible to have a (somewhat irregular) statistical
model for
which consistency holds even though the Fisher information does not go
to infinity.

\subsection{Comparison of rates}

Now, we can compare these optimal rates with the rates obtained in Theorem
\ref{TD2}. Doing as such, we compare a semiparametric model with a parametric
sub-model. However, a minimax rate for a given parameter in a semiparametric
model cannot be faster than the rate obtained for any parametric sub-model,
hence the previous results are bounds for the rates in the general model
considered in this paper.

Neglecting the logarithmic terms, and considering only the estimation
of~$\beta_{i}$ for $i=1,2$, the rates above are $\Delta_{n}^{\gamma_{i}}$,
whereas in Theorem~\ref{TD2}, and upon choosing $u_{n}$ optimally [i.e.,
$\rho$ as large as possible in~(\ref{D31})], they are $\Delta_{n}%
^{\gamma^{\prime}_{i}}$, where
\[
\gamma_{i}=\frac{2\beta_{i}-\beta_{1}}4,\qquad \gamma^{\prime}_{i}=\cases{
\displaystyle \gamma_{i} \frac{2}{2+\beta_{1}}-\varepsilon, &\quad if $\beta_{1}%
\leq\bigl(\sqrt{97}-1\bigr)/6\approx1.475$,\vspace*{2pt}\cr
\displaystyle \gamma_{i} \frac{8}{5\beta_{1}+3\beta_{1}^{2}}-\varepsilon, &\quad if
$\beta_{1}\geq\bigl(\sqrt{97}-1\bigr)/6$,}
\]
and $\varepsilon>0$ arbitrarily small (and if $\beta_{i}>\beta_{1}/2$ when
$i=2$).

As it should be, we have $\gamma_{i}\leq\gamma_{i}^{\prime}$, and if equality
were holding we would conclude that our estimators achieve the minimax rate
(up to $\Delta_{n}^{-\varepsilon}$, of course, but $\varepsilon$ is
arbitrarily small). What one can say is that the actual minimax rate lies
somewhere in between these two values, and the ratio $\gamma_{i}/\gamma
_{i}^{\prime}$ is a kind of (imperfect) measure of the quality of the
estimators proposed in Section~\ref{secD2}: the closest to $1$, the closest
to optimality. Then we can conclude the following:

(a) This ratio is the same for $j=1,2$, which is an a priori
surprising result: the quality of our estimator for $\beta_{2}$,
relative to
the optimal estimators in the stable sub-model, is the same as for
$\beta_{1}$.

(b) This ratio is close to $1$ (near optimality) when $\beta_{1}$ is small,
and decreases down to $4/11$ as $\beta_{1}$ increases up to $2$. The worst
value is small, but not catastrophically such, especially in light of
the fact
that we are considering semiparametric estimators whereas the rates are
optimal in the parametric context (i.e., assuming additional structure).

\section{Simulation results}\label{secMC}

We now provide some simulation evidence regarding the estimators in the case
where $j=2$; we are attempting to estimate the first two jump activity indices
of the process $\beta_{1}$ and $\beta_{2}$. The data generating process
is a
stochastic volatility model for $X_{t}$ with jumps driven by two stable
processes $Y^{1}$ and $Y^{2}$, with $W,Y^{1},Y^{2}$ independent below:
%
\begin{equation}\label{eqMCsde}%
dX_{t}=\sigma_{t}\,dW_{t}+\theta_{1}\,dY_{t}^{1}+\theta_{2}\,dY_{t}^{2}
\end{equation}
with $\sigma_{t}=v_{t}^{1/2}$, $dv_{t}=\kappa(\eta-v_{t})\,dt+\gamma v_{t}
^{1/2}\,dB_{t}+dJ_{t}$, $\mathbb{E}[dW_{t}\,dB_{t}]=\rho \,dt$, $\eta^{1/2}=0.25$,
$\gamma=0.5$, $\kappa=5$, $\rho=-0.5$, the volatility jump term $J$ is a
compound Poisson jump process with jumps that are uniformly distributed on
$[-0.3,0.3]$ and intensity $\lambda=10$ and $X_{0}=1$. Recall that the second
component can be identified only if $\beta_{2}>\beta_{1}/2$. We
consider the
situation where $(\beta_{1},\beta_{2})=(1.00,0.75)$.

Given $\eta$, each scale parameter $\theta_{i}$ (or equivalently $A_{T}^{i}$)
of the stable process in simulations is calibrated to deliver different
various values of the tail probability $P_{i}=\mathbb{P}(|\Delta Y_{t}%
^{i}|\geq4\eta^{1/2}\Delta_{n}^{1/2})$. In the various simulations'
design, we
hold $\eta$ fixed and consider the cases where $P_{1}=0.05$ and $P_{2}=0.005$.
We sample the process $X$ over $T=21$ days ($6.5$ hours per day) every
$\Delta_{n}=0.01$ second. This results of course in a number of observations
(nearly $5\times10^{7}$) that is unrealistically high for most high-frequency
financial data series, at least presently, but extremely large numbers of
observations are needed if we are going to be able to see the component
$\beta_{2}$ of the model ``behind''
the two
components with indices of activity $2$ (the continuous component) and
$\beta_{1}$ (the most active jump component). Of course, much smaller datasets
would be sufficient in the absence of a continuous component.%

Note that in general, and besides the preliminary estimators $\widetilde
{\beta}^{n}_{i}$ and $\widetilde{\Gamma}^{n}_{i}$, we need to choose the
number $\alpha>0$ coming in the definition of the set $D_{n}$. Since in
practice $n$ (or $\Delta_{n}$) is given, we need to choose in fact the number
$\alpha u_{n}^{\eta}$. So in concrete situations one probably can
forget about
the preliminary estimators and take a domain $D_{n}$ which is the set
of all
$(x_{i},\gamma_{i})$ in $D$ with $\gamma_{i}\leq A$ for some ``reasonably
chosen'' $A$, or even $A=\infty$.

This is what we do below, by taking the estimators to be
%
\begin{eqnarray}\label{eqargmin4}
&&
(\overline{\beta}{}^{\prime n}_{1},\overline{\beta}{}^{\prime n}_{2}%
,\overline{\Gamma}{}^{\prime n}_{1},\overline{\Gamma}{}^{\prime n}_{2}%
)\nonumber\\[-8pt]\\[-8pt]
&&\qquad=\mathop{\arg\min}_{( x_{1},\gamma_{1},x_{2},\gamma_{2})
}\sum_{l=1}^{L}\biggl( U(v_{l}u_{n},\Delta_{n})_{T}-\frac{\gamma_{1}}%
{(v_{l}u_{n})^{x_{1}}}-\frac{\gamma_{2}}{(v_{l}u_{n})^{x_{2}}}\biggr)
^{2},\nonumber
\end{eqnarray}
where the cutoff levels $v_{l}u_{n}$ are chosen in terms of the number
$\alpha_{l}$ of the long-term standard deviation $\sqrt{\eta\Delta_{n}}$
over a time lag $\Delta_{n}$ of the continuous martingale part of the process:
we take $\alpha_{l}$ to be $\{ 7,10,15,20\} $ and multiples
$\{ 2,4,6\} $ thereof (giving all together $L=10$ distinct
values). Here we know $\eta$: we could also estimate for each path the average
volatility, using truncated estimators for the integrated volatility [see,
e.g., \citet{mancini04} and \citet{yacjacod09a}].

The optimization problem~(\ref{eqargmin4}) is a quadratic problem
similar to
classical nonlinear least squares minimization. In situations where the
parameter space is high dimensional, the objective function can exhibit local
extrema, which can make the search for the optimal solution
time-consuming as
many starting values must be employed to validate the solution. In the
case of
the application here, we are only including $4$ parameters, and for
this small
dimension, this is not causing many difficulties. In any case, it is unlikely,
given the slow rates of convergence, that one would want to go beyond the
second index $\beta_{2}$ in practice.

%
\begin{figure}

\includegraphics{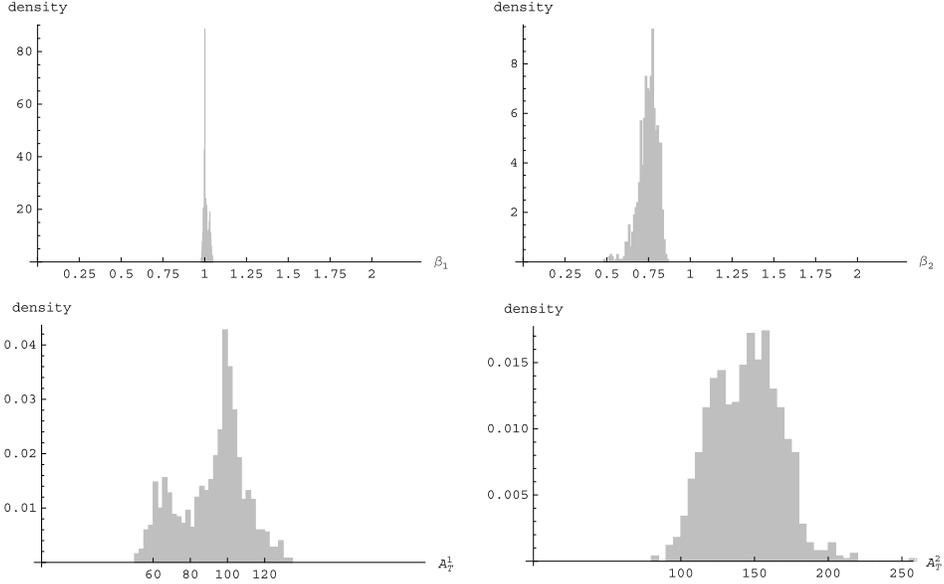}

\caption{Monte Carlo simulation results: estimators $\widehat{\beta}_{1}^{n}$
(upper left graph), $\widehat{\beta}_{2}^{n}$ (upper right graph), $\widehat{A}%
_{T}^{1,n}$ (lower left graph), $\widehat{A}_{T}^{2,n}$ (lower right
graph).}%
\label{figHistBeta1and2}
\end{figure}

The results in Figure~\ref{figHistBeta1and2} are obtained with $M=1000$
simulations: the estimators appear to be reasonably good, but then
again this
is for an unrealistically large number of observations, at least from the
point of view of financial applications; it is perhaps feasible in other
applications, such as Internet data traffic or wind measurement.

\section{Conclusions}\label{conclusion}

This paper determined theoretically what the successive BG indices
are and how they are identified, including the perhaps surprising
theoretical bound on the
identification of the successive indices as a function of the previous ones.
This result clarifies the border between the aspects of the jump
measure which
are identifiable from those which are not on the basis of discrete
observations on a finite time horizon. Beyond the leading index, the
identification requires in practice vast quantities of data which are
out of
reach of financial applications at present but may be relevant in other fields
(such as the study of turbulence data, or Internet traffic). We showed through
explicit calculations of Fisher's information that this limitation is a
genuine, inescapable feature of the problem. There are a number of important
questions that this paper does not touch upon: central limit theorems
for the
estimators, estimators that achieve the optimal rates of convergence,
estimators that are robust to microstructure noise, estimators that are
applicable with random sampling intervals, among others. The issue of the
optimality of the rates in general remains an open question.

\begin{appendix}\label{app}

\section*{Appendix: Proofs}

We use the following notation throughout the \hyperref[app]{Appendix}. First, $K$
denotes a
constant which may change from line to line, and may depend on the
characteristics or the law of the processes at hand. It never depends
on $n$,
and it is denoted as $K_{p}$ if it depends on an additional parameter $p$.
Second, for any sequence $Z_{n}$ of variables and any sequence $v_{n}$ of
positive numbers,
%
\begin{equation}\label{S8}
Z_{n}=\cases{
O_{P}(v_{n}), &\quad if $Z_{n}/v_{n}$ is bounded in probability,\cr
o_{P}(v_{n}), &\quad if $Z_{n}/v_{n}\stackrel{\mathbb{P}}{\longrightarrow
}0$.}\
\end{equation}

\setcounter{section}{0}
\section{\texorpdfstring{Proof of Theorem \protect\ref{TLC1}}{Proof of Theorem 1}}

(1) We fix\vspace*{1pt} $F\in\mathcal{T}_{3}^{(1)}$, with $F$ given by~(\ref{400}).
We also
consider another $F^{\prime}\in\mathcal{T}_{3}^{(1)}$, with $F^{\prime
}$ given
by~(\ref{400}) with $\beta_{i}^{\prime}$, $a_{i}^{\prime}$ and
$\widetilde
{F}^{\prime}$. As said before, it is not a restriction to assume the
representation~(\ref{400}) with the same $\eta>0$ for both~$F$ and
$F^{\prime
}$. Set
%
\begin{equation}\label{405}%
j=\inf\bigl(1\leq1\dvtx (\beta_{i},a_{i})\neq(\beta_{i}^{\prime},a_{i}^{\prime})\bigr).
\end{equation}
The result amounts to proving the following two properties, with $j$ as above
and $b,b^{\prime}\in\mathbb{R}$ and $c,c^{\prime}\geq0$:
%
\begin{eqnarray}
\label{402}%
\beta_{j}&\geq&\frac{\beta_{1}}{2} \quad\Rightarrow\quad Q_{b,c,F}\perp
Q_{b^{\prime},c^{\prime},F^{\prime}},\hspace*{-35pt}
\\
\label{403}%
\beta_{j}&<&\frac{\beta_{1}}{2}\quad\Rightarrow\quad\cases{
\exists b^{\prime\prime}\in\mathbb{R},\exists F^{\prime\prime}\in
\mathcal{T}_{3}^{(1)}\cr
\qquad\mbox{with }F^{\prime\prime}=F^{\prime}\mbox{ on }[-\eta,\eta)\mbox{ and }%
Q_{b,c,F}\not\perp Q_{b^{\prime\prime},c,F^{\prime}}.}\hspace*{-35pt}
\end{eqnarray}
These conditions being symmetrical in $F$ and $F^{\prime}$, in both
(\ref{402}) and~(\ref{403}) we may assume
%
\begin{equation} \label{404}%
\mbox{either }\beta_{j}>\beta_{j}^{\prime}\quad\mbox{or}\quad\beta_{j}=\beta
_{j}^{\prime}\quad\mbox{and}\quad a_{j}>a_{j}^{\prime}.
\end{equation}

(2) In this step we assume~(\ref{404}). We set
\[
\widehat{F}(dx)=\sum_{i\geq1}\frac{a_{i}\beta_{i}}{|x|^{1+\beta_{i}}
} 1_{[-\eta,\eta]}(x)\,dx,\qquad \widehat{F}^{\prime}(dx)= \sum_{i\geq1}%
\frac{a_{i}^{\prime}\beta_{i}^{\prime}} {|x|^{1+\beta_{i}^{\prime}}%
} 1_{[-\eta,\eta]}(x)\,dx.
\]
Then $\widehat{F}^{\prime}=f\bullet\widehat{F}$, where $f=\frac
{g^{\prime}}%
{g}$ (with $\frac00=1$) and $g=H+G$ and $g^{\prime}=H+G^{\prime}$ and
\begin{eqnarray*}
H(x)&=&\sum_{i=1}^{j-1}\frac{a_{i}\beta_{i}}{|x|^{\beta_{i}}} 1_{[-\eta
,\eta
]}(x),\qquad G(x)=\sum_{i\geq j}\frac{a_{i}\beta_{i}}{|x|^{\beta_{i}} }%
1_{[-\eta,\eta]}(x),\\
G^{\prime}(x)&=&\sum_{i\geq j}\frac{a_{i}^{\prime}%
\beta_{i}^{\prime}}{|x|^{\beta_{i}^{\prime}}} 1_{[-\eta,\eta]}(x).
\end{eqnarray*}
On $[-\eta,\eta]$ we have $f-1=\frac{G^{\prime}-G}{H+G}$ and
\begin{eqnarray*}
&&
G(x)-G^{\prime}(x)\\
&&\qquad=\frac{a_{j}\beta_{j}}{|x|^{1+\beta_{j}}}\biggl(
1-\frac{a_{j}^{\prime}\beta_{j}^{\prime}}{a_{j}\beta_{j}} |x|^{\beta_{j}
-\beta_{j}^{\prime}}+\sum_{i\geq j+1}\frac{a_{i}\beta_{i}}{a_{j}\beta_{j}
} x^{\beta_{j}-\beta_{i}}-\sum_{i\geq j+1}\frac{a_{i}^{\prime} \beta
_{i}^{\prime}}{a_{j}\beta_{j}} |x|^{\beta_{j}-\beta_{i}^{\prime}}\biggr) .
\end{eqnarray*}
By virtue of (ii), (iii) and (iv) of~(\ref{400}), and of~(\ref{404}),
we then
deduce that
%
\begin{equation}\label{406}%
x\in(-\varepsilon,\varepsilon) \quad\Rightarrow\quad\cases{
\displaystyle A_{-} |x|^{\beta_{1}-\beta_{j}}\leq|f(x)-1|\leq A_{+} |x|^{\beta_{1}
-\beta_{j}},\vspace*{2pt}\cr
\displaystyle \frac{A_{-}}{|x|^{1+\beta_{1}}} \leq g(x) \leq
\frac{A_{+}}{|x|^{1+\beta_{1}}},}
\end{equation}
for three constants $A_{+}>A_{-}>0$ and $\varepsilon\in(0,\eta)$,
depending on
the two sequences $(\beta_{i},a_{i})$ and $(\beta_{i}^{\prime
},a_{i}^{\prime
})$.

(3) Now we\vspace*{1pt} prove~(\ref{402}). Since
$\alpha(F,F^{\prime})\geq\alpha (\widehat {F}, \widehat{F}^{\prime})$,
it is enough to show that $\alpha(\widehat {F},\widehat{F}^{\prime})
=\infty$. By~(\ref{406}), $|f(x)-1|\leq1$ when
$x\in(-\varepsilon^{\prime},\varepsilon^{\prime})$ for some
$\varepsilon ^{\prime}\in(0,\varepsilon]$. Thus
\[
\alpha(\widehat{F},\widehat{F}^{\prime}) \geq \int_{-\varepsilon
}^{\varepsilon^{\prime} }|f(x)-1|^{2} g(x)\,dx \geq A_{-}^{3} \int
_{-\varepsilon^{\prime}}^{\varepsilon^{\prime}} |x|^{\beta_{1}-2\beta
_{j}%
-1}\,dx.
\]
The last integral is infinite when $\beta_{j}\geq\beta_{1}/2$, and~(\ref{403})
follows by~(\ref{LC1}).\vspace*{1pt}

(4) Finally we prove~(\ref{403}). Recall that $F=\widehat{F}+\widetilde
{F}$ and
$F^{\prime}=\widehat{F}^{\prime}+\widetilde{F}^{\prime}$. The measure
$F^{\prime\prime}=\widehat{F}^{\prime}+\widetilde{F}$ is obviously in
$\mathcal{T}^{(1)}_{3}$ and satisfies $F^{\prime\prime}=f\bullet F$. Since
$f(x)=1$ outside $[-\eta,\eta]$, the quantity $\alpha^{\prime
}(F,F^{\prime
\prime})$ introduced in~(\ref{LC1}) is
\begin{eqnarray*}
\alpha^{\prime}(F,F^{\prime\prime}) & = & \int_{-\eta}^{\eta}%
x\bigl(f(x)-1\bigr)g(x)\,dx\\
& \leq & A_{+}^{2}\int_{-\varepsilon}^{\varepsilon}|x|^{-\beta_{j}}\,dx
+\biggl(\int_{\varepsilon}^{\eta}+\int_{-\eta}^{-\varepsilon}%
\biggr) |x| |f(x)-1| g(x)\,dx,
\end{eqnarray*}
which is finite by~(\ref{406}) (because $\beta_{j}<\beta_{1}/2<1$) and because
$f$ and $g$ are bounded on $[\varepsilon,\eta]\cup[-\eta,-\varepsilon]$.
Therefore the number $b^{\prime\prime}=b-\int_{0}^{\eta\wedge1}%
x(f(x)-1)g(x)\,dx$ is well defined. Now we consider the two triples $(b,c,F)$
and $(b^{\prime\prime},c,F^{\prime\prime})$. From what precedes they satisfy
the first and the last three properties in~(\ref{LC1}). We also have
by~(\ref{406})
\begin{eqnarray*}
\alpha(F,F^{\prime\prime}) & = & \int_{-\eta}^{\eta}\bigl(|f(x)-1|^{2}%
\wedge|f(x)-1|\bigr) g(x)\,dx\\
& \leq & A_{+}^{2}\int_{-\varepsilon}^{\varepsilon} \bigl( |x|^{-\beta
_{j}-1}\wedge( A_{+}|x|^{\beta_{1}-2\beta_{j}-1}) \bigr) \,dx\\
&&{} +\biggl(\int_{\varepsilon}^{\eta}+\int_{-\eta}^{-\varepsilon
}\biggr) \bigl(|f(x)-1|^{2}\wedge|f(x)-1|\bigr) g(x)\,dx.
\end{eqnarray*}
Since $\beta_{j}<\beta_{1}/2$ and that $f$ and $g$ are bounded on
$[\varepsilon,\eta]\cup[-\eta,-\varepsilon]$, we deduce $\alpha
(F,F^{\prime
\prime})<\infty$. So all conditions in~(\ref{LC1}) are satisfied, and
we have
proved~(\ref{403}).

\section{Comparing big jumps and big increments}

Before starting, let us mention that for the proofs of Theorems \ref
{TD1} and
\ref{TD2} one may use a localization argument which allows us to replace
Assumption~\ref{A2} by the so-called ``strengthened Assumption~\ref{A2},''
which is the same except that all processes $b_{t}$, $c_{t}$,
$a^{i}_{t}$ are
bounded, as well as the process $A^{j+1}_{t}$ and $X_{t}$ itself.

In this section we compare the number of ``large'' increments of $X$
with the
number of correspondingly large jumps, that is, the numbers
%
\begin{equation}
\label{S3}
V(u)_{t} = \sum_{s\leq t} 1_{\{|\Delta X_{s}|>u\}}.%
\end{equation}
We will indeed show that the difference $U(u_{n},\Delta_{n})_{T}-
V(u_{n}%
)_{T}$ is negligible for our purposes, when the sequence $u_{n}$ satisfies
(\ref{D31}). The reason for doing this is that the analysis of the processes
$V(u_{n})$ is an easy task. Indeed, as soon as $u_{n}\to0$,
%
\begin{equation}
\label{KP}
V(u_{n})_{T}-\overline{A}(u_{n})_{T} = O_{P}(u_{n}^{-\beta_{1}/2}).
\end{equation}
To see this, we observe that each process $M^{n}=u_{n}^{\beta_{1}/2}(
V(u_{n})-\overline{A}(u_{n})) $ is a~quasi-left continuous, purely
discontinuous, martingale with jumps smaller than $u_{n}^{\beta
_{1}/2}$, which
goes to $0$. Its predictable quadratic variation is $\langle M^{n}
,M^{n}\rangle=u_{n}^{\beta_{1}} \overline{A}(u_{n})$, which by~(\ref{31})
converges for each $t$ to $A_{t}^{1}$. Since further $A^{1}$ is a~continuous
process, it follows from Theorem VI.4.13 of \citet{jacodshiryaev2003}, for
example, that the sequence $M^{n}$ is $C$-tight (and even converges in
law), so
a fortiori,~(\ref{KP})~holds.

The main result of this section is the next proposition:
%
\begin{proposition}
\label{PD1} Under the strengthened Assumption~\ref{A2} and if the sequence
$u_{n}$ satisfies~(\ref{D31}), we have
%
\begin{equation} \label{D30}%
U(u_{n},\Delta_{n})_{T}-V(u_{n})_{T} = \frac{1}{u_{n}^{\beta_{1}}}
O_{P}(u_{n}^{\beta_{1}-\beta_{j+1}}+u_{n}^{\beta_{1}/2} ).
\end{equation}
\end{proposition}

The proof is based on a series of lemmas. The constant $K$ may depend
on an
implicit way on the bounds in this strengthened assumption, but not on
the two
numbers $u,r\in(0,1)$ which are fixed in most of this section.

With any c\`{a}dl\`{a}g process $Y$ and $u\in(0,1]$, we associate the process
and the variables
%
\begin{equation} \label{D21}%
Y(u)_{t} = \sum_{s\leq t}\Delta Y_{s}1_{\{|\Delta Y_{s}|>u\}},\qquad
\zeta(Y,u)_{i}^{n} = 1_{\{|\Delta_{i}^{n}Y|>u\}}.
\end{equation}
For simpler notation, we denote by $\mathbb{E}_{i-1}^{n}$ and $\mathbb{P}
_{i-1}^{n}$, respectively, the conditional expectation and conditional
probability, with respect to $\mathcal{F}_{(i-1)\Delta_{n}}$.

\begin{lemma}
\label{LD1} For all $u,r\in(0,1]$ with $u^{r}<1/3$, all $w\in(0,1/3)$
and all
$k\geq1$, we have
%
\begin{eqnarray} \label{D23}%
&&\mathbb{P}_{i-1}^{n}\bigl(\Delta_{i}^{n}V(u)\geq k\bigr) \leq( K\Delta
_{n} u^{-\beta_{1}}) ^{k},\hspace*{-35pt}
\\
%
\label{D25}
&&\mathbb{P}_{i-1}^{n}\bigl(u(1-w)<\Delta_{i}^{n}X(u^{1+r})\leq
u(1+w)\bigr)\hspace*{-35pt}\nonumber\\[-8pt]\\[-8pt]
&&\qquad\leq K\bigl( \Delta_{n}u^{-\beta_{1}}w+\Delta_{n}u^{- \beta_{j+1}}
+\Delta_{n}^{2}u^{-\beta_{1}(2+r)}+\Delta_{n}^{3}u^{-\beta_{1}(3+3r)}\bigr).
\nonumber\hspace*{-35pt}
\end{eqnarray}
Moreover there is a $\gamma>0$ such that, if
%
\begin{equation} \label{D265}
\Delta_{n} \leq\gamma u^{\beta_{1}(1+r)},%
\end{equation}
we have for all $u\in(0,1]$
%
\begin{equation} \label{D24}
\mathbb{E}_{i-1}^{n}\bigl( \vert\zeta(X(u^{1+r}),u)_{i}^{n}-\Delta
_{i}^{n}V(u)\vert\bigr) \leq K\bigl( \Delta_{n}^{2} u^{-\beta
_{1}(2+r)}+\Delta_{n}^{3} u^{-\beta_{1}(3+3r)}\bigr) .\hspace*{-35pt}
\end{equation}
\end{lemma}
\begin{pf}
If $D\subset\mathbb{R}$ the compensator of the process $N(D)_{t}= \sum
_{s\leq
t}1_{D}(\Delta X_{s})$ is $\int_{0}^{t}F_{s}(D)\,ds$. Our strengthened
assumption implies the existence of a constant $\theta$ such that
$F_{s}(D)\leq\phi(D)$, where
\[
\phi(D)=\cases{
\theta u^{-\beta_{1}}, &\quad if $D\subset[-u,u]^{c}$,\vspace*{1pt}\cr
\theta( u^{-\beta_{1}}w+u^{-\beta_{j+1}}), &\quad if
$D=\bigl[-u(1+w),-u\bigr)\cup\bigl(u,u(1+w)\bigr]$,\vspace*{1pt}\cr
&\quad $0<w\leq1$.}
\]
Then for any finite stopping time $S$ we have
\[
\mathbb{E}\bigl(N(D)_{S+t}-N(D)_{S}\mid\mathcal{F}_{S}\bigr) \leq t \phi(D).
\]
Let $S(D)_{0}=(i-1)\Delta_{n}$ and $S(D)_{1},S(D)_{2},\ldots$ be the
successive jump times of $N(D)$ after time $(i-1)\Delta_{n}$. What precedes
implies that for $k\geq1$ and on the set $\{S(D)_{k-1}<i\Delta_{n}\}$,
\[
\mathbb{P}\bigl(S(D)_{j}\!\leq\!i\Delta_{n}\mid\mathcal{F}_{S(D)_{k-1}}\bigr)\!\leq\!
\mathbb{E}\bigl( N(D)_{i\Delta_{n}}\!-\!N(D)_{S(D)_{k-1}}\mid\mathcal{F}
_{S(D)_{k-1}}\bigr)\!\leq\!\Delta_{n} \phi(D).
\]
An induction on $k$ yields the following, which gives us the first part of
(\ref{D23}):
%
\begin{equation}\label{D250}%
\mathbb{P}_{i-1}^{n}\bigl(\Delta_{i}^{n}N(D)\geq k\bigr) = \mathbb{P}_{i-1}^{n}
\bigl(S(D)_{k}\leq i\Delta_{n}\bigr) \leq( \Delta_{n} \gamma(D)) ^{k}.
\end{equation}

In the same way, if $D\cap D^{\prime}=\varnothing$, the set $\{\Delta_{i}
^{n}N(D)\geq k,\Delta_{i}^{n}N(D^{\prime})\geq1\}$ is the union for
$l=1,\ldots,k+1$ of the sets $\Gamma_{l}=\{S(D)_{l-1}<S(D^{\prime}
)_{1}<S(D)_{l}\leq i\Delta_{n}\}$, whereas
\begin{eqnarray*}
\hspace*{-4pt}&&
\mathbb{P}_{i-1}^{n}\bigl(S(D)_{l-1}<S(D^{\prime} )_{1}<S(D)_{l}\leq i\Delta
_{n}\bigr)\\
\hspace*{-4pt}&&\quad=\mathbb{E}_{i-1}^{n}\bigl( 1_{S(D)_{l-1} <S(D^{\prime})_{1}<i\Delta
_{n}}\mathbb{P}\bigl( N(D)_{i\Delta_{n} } - N(D)_{S(D^{\prime})_{1}%
}\geq k-l+1 \hspace*{-0.3pt}\mid\hspace*{-0.3pt}\mathcal{F}_{S(D)_{1}} \bigr) \bigr) \\
\hspace*{-4pt}&&\quad\leq(\Delta_{n}\phi(D))^{k-l+1}\mathbb
{P}_{i-1}^{n}\bigl(S(D)_{l-1}<S(D^{\prime
})_{1}<i\Delta_{n}\bigr)\\
\hspace*{-4pt}&&\quad=(\Delta_{n}\phi(D))^{k-l+1}\\
\hspace*{-4pt}&&\qquad{}\times\mathbb{E}_{i-1}^{n}\bigl( 1_{S(D)_{l-1}
<i\Delta_{n}}\mathbb{P}\bigl( N(D^{\prime})_{i\Delta_{n}}- N(D^{\prime
})_{S(D^{\prime\prime})_{l-1}}\geq1 \mid\mathcal{F}_{S(D)_{1}} \bigr)
\bigr) \\
\hspace*{-4pt}&&\quad\leq(\Delta_{n}\phi(D))^{k-l+1}\Delta_{n}\phi(D^{\prime})\mathbb{P}
_{i-1}^{n}\bigl(S(D)_{l-1}<i\Delta_{n}\bigr),
\end{eqnarray*}
where~(\ref{D250}) has been applied twice. Another application of the same
then yields
%
\begin{eqnarray} \label{D231}%
D\cap D^{\prime}&=&\varnothing\Rightarrow\mathbb{P}_{i-1}^{n}\bigl(\Delta_{i}
^{n}N(D)\geq k,
\Delta_{i}^{n}N(D^{\prime})\geq1\bigr)\nonumber\\[-8pt]\\[-8pt]
&\leq&(k+1)\Delta_{n}
^{k+1}\phi(D)^{k}\phi(D^{\prime}).\nonumber
\end{eqnarray}

Next, let $w\in(0,1/3]$. By convention $(a,b]=\varnothing$ when $a\geq
b$ below.
If $u(1-w)<\Delta_{i}^{n} X(u^{1+r})\leq u(1+w)$ we have four (nonexclusive)
possibilities: either $\Delta_{i}^{n}N((u^{1+r},\infty))\geq3$, or
$\Delta
_{i}^{n}N((u^{1+r},u/3])=\Delta_{i} ^{n}N((u/3,\infty))=1$, or $\Delta_{i}
^{n}N((u/3,\infty))=2$, or $\Delta_{i}^{n}N((u(1-w),u(1+w)])=1$. We an
analogous implication if $-u(1+w)<\Delta_{i}^{n} X(u^{1+r})\leq
-u(1-w)$. Then
(\ref{D25}) easily follows from~(\ref{D250}) applied with $D=[-u^{1+r}%
,u^{1+r}]^{c}$, with $D=[-u/3,u/3]^{c}$ and with $D=
[-u(1+w),-u(1-w))\cup
(u(1-w),u(1+w)]$, and from~(\ref{D231}) applied with $D=(-u/3,-u^{1+r}%
)\cup(u^{1+r},u/3]$ and $D^{\prime}=[-u/3,u/3]^{c}$.

Finally we prove~(\ref{D24}). Let $H=|\zeta(X(u^{1+r}),u)_{i}^{n}-\Delta
_{i}^{n}V(u)|$ and $D=[-u/2$, $-u^{1+r})\cup(u^{1+r},u/2]$ and $D^{\prime
}=[-u/2,u/2]^{c}$ and $D^{\prime\prime}=D\cup D^{\prime}$. From what precedes,
we have
%
\begin{eqnarray}\label{D240}
\mathbb{P}_{i-1}^{n}\bigl(\Delta_{i}^{n}N(D^{\prime\prime})\geq k\bigr)&\leq&\bigl(
\theta\Delta_{n}u^{-\beta_{1}(1+r)}\bigr) ^{k},\nonumber\\
\mathbb{P}_{i-1}^{n}\bigl(\Delta_{i}^{n}N(D^{\prime})=2\bigr)&\leq&\theta^{2}\Delta
_{n}^{2}u^{-2\beta_{1}},\\
\mathbb{P}_{i-1}^{n}\bigl(\Delta_{i}^{n}N(D)=\Delta_{i}^{n}N(D^{\prime}
)=1\bigr)&\leq&\theta^{2}\Delta_{n}^{2}u^{-\beta_{1}(2+r)}.
\nonumber
\end{eqnarray}
We have $H=0$ on the sets $\{\Delta_{i}^{n}N(D^{\prime\prime})\leq1\}$ and
$\{\Delta_{i}^{n}N(D^{\prime\prime})=\Delta_{i}^{n}N(D)=2\}$, and
$H\leq k-1$
on the set $\{\Delta_{i}^{n}N(D^{\prime\prime})=k\}$, for all $k\geq2$. Thus
if $v=\theta\Delta_{n}u^{-\beta_{1}(1+r)}$,
\begin{eqnarray*}
\mathbb{E}_{i-1}^{n}(H) & \leq & \sum_{k=3}^{\infty}k \mathbb{P}_{i-1}
^{n}\bigl(\Delta_{i}^{n}N(D^{\prime\prime})\geq k\bigr)+\mathbb
{P}_{i-1}^{n}\bigl(\Delta
_{i}^{n}N(D^{\prime})=2\bigr)\\
&&{} +\mathbb{P}_{i-1}^{n}\bigl(\Delta_{i}^{n}N(D)=\Delta_{i}^{n}N(D^{\prime
})=1\bigr)\\
& \leq & \sum_{k=3}^{\infty}kv^{k}+\theta^{2}\Delta_{n}^{2}u^{-2\beta_{1}
}+\theta^{2}\Delta_{n}^{2}u^{-\beta_{1}(2+r)}%
\end{eqnarray*}
by~(\ref{D240}). When $v\leq1/2$, that is, when $\Delta_{n}\leq\gamma
u^{\beta_{1}(1+r)}$ for $\gamma=1/2\theta$, we have $\sum_{k=3}^{\infty}
kv^{k}\leq Kv^{3}$, and the above is smaller than the right-hand side
of~(\ref{D24}).
\end{pf}

\begin{lemma}
\label{LD301} Let $q\geq2$ and $u,r\in(0,1)$. As soon as~(\ref{D265}) holds
for some constant $\gamma>0$, we have
%
\begin{equation} \label{D266}%
\mathbb{E}\bigl(\bigl|\Delta_{i}^{n}\bigl(X-X(u^{1+r})\bigr)\bigr|^{q}\bigr) \leq K_{\gamma,q}
\Delta_{n}\bigl(\Delta_{n}^{q/2-1}+u^{(q-\beta_{1})(1+r)}\bigr).
\end{equation}
\end{lemma}
\begin{pf}
Letting $X^{c}$ and $\mu$ be the continuous martingale part and the jump
measure of $X$, we have $X-X(u^{1+r})=B+B^{\prime}+X^{c}+M$, where
\[
B_{t}^{\prime}=-\int_{0}^{t}ds\int_{\{u^{1+r}<|x|\leq1\}}xF_{s}(dx),\quad
M_{t}=\int_{0}^{t}\int_{\{0<|x|\leq u^{1+r}\}}x (\mu-\nu)(ds,dx).
\]
By the strengthened Assumption~\ref{A2}, for any $y>0$ the integral
$\int_{\{|x|>y\}}|x| F(dx)$ is smaller than $K$ when $\beta_{1}<1$, than
$K\log\frac{1}{y}$ when $\beta_{1}=1$, and than $Ky^{1-\beta_{1}}$ when
$\beta_{1}>1$. Therefore, since~(\ref{D265}) implies $2\beta_{1}%
(1+r)>(\beta_{1}-1)^{+}$ we have $|\Delta_{i}^{n}B^{\prime}|\leq
K_{\gamma
}\sqrt{\Delta_{n}}$. The strengthened Assumption~\ref{A2} also implies
$|\Delta^{n}_{i}B|\leq K\Delta_{n}$ and, by well-known estimates about
continuous and purely discontinuous martingales [see, e.g.,
\citet{yacjacod11}], we also deduce that
\[
\mathbb{E}(|\Delta_{i}^{n}M^{\prime}|^{q})\leq K_{q}\Delta_{n}
u^{(q-\beta
_{1})(1+r)} ,\qquad \mathbb{E}(|\Delta_{i}^{n}X^{c}|^{q})\leq K_{q}\Delta
_{n}^{q/2}.
\]
All these estimates readily give~(\ref{D266}).
\end{pf}
\begin{pf*}{Proof of Proposition~\ref{PD1}}
(a) It follows from~(\ref{D31}) that
$u_{n}^{2\beta_{1}}/\Delta_{n} \to\infty$, so for any $\gamma>0$~(\ref{D265})
is satisfied for all $r\in(0,1)$ and all $n$ large enough. Hence, both
estimates~(\ref{D24}) and~(\ref{D266}) hold, with constants $K$ and
$K_{\gamma,q}$ independent of $r$, for all $n$ large enough.

The following inequality, where $u,w\in(0,1)$ and $x,y\in\mathbb{R}$, is
elementary:
\[
\bigl\vert1_{\{x+y>u\}}-1_{\{x>u\}}\bigr\vert \leq1_{\{|y|\geq
uw\}}+1_{\{u(1-w)<|x|\leq u(1+w)\}}.
\]
We apply this with $x=\Delta_{i}^{n}X(u^{1+r})$ and $x+y=\Delta
^{n}_{i}X$ and
$u=u_{n}$, and with $w\leq1/3$ to be chosen later.\vadjust{\goodbreak} In order to evaluate the
probabilities for having $|y|\geq u_{n}w$, respectively,
$u_{n}(1-w)<|x|\leq
u_{n}(1+w)$, we use~(\ref{D266}) and Markov's inequality, respectively,
(\ref{D25}).
This gives that $\mathbb{E} (|\zeta(X(u_{n}^{1+r}),u_{n})^{n}_{i}-
\zeta(X,u_{n})^{n}_{i}|)$ is smaller, for all $q\geq2$, than
\[
\frac{K_{q}\Delta_{n}}{u_{n}^{\beta_{1}}}\biggl(
\frac{\Delta_{n}^{q/2-1}}{w^{q} u_{n}^{q-\beta_{1}}}+ \frac{u_{n}
^{(q-\beta_{1})r}}{w^{q}}+w+\frac{\Delta_{n}} {u_{n}^{\beta_{1}(1+r)}}%
+\frac{\Delta_{n}^{2}}{u_{n}^{\beta_{1}(2+3r)}} +u_{n}^{\beta_{1}-\beta
_{j+1}%
}\biggr).
\]
Optimizing over $w$ leads to take $w=w_{n}$ such that $w_{n}^{q+1}=
u_{n}^{(q-\beta_{1})r}+\Delta_{n}^{q/2-1}/\allowbreak u_{n}^{q-\beta_{1}}$, which is
indeed smaller than $1/3$ for all $n$ large. Thus, putting the above together
with~(\ref{D24}), and recalling that $\Delta_{n}\leq Ku_{n}^{1/\rho}$,
we end
up with
%
\begin{equation}\label{D119}%
\mathbb{E}\bigl(|\zeta(X,u_{n})^{n}_{i}-\Delta^{n}_{i}V(u_{n})|\bigr) \leq \frac
{K_{q}\Delta_{n}}{u_{n}^{\beta_{1}}}\sum_{k=1}^{5}u_{n}^{x_{k}},
\end{equation}
where $x_{k}=x_{k}(q,r)$ are given by
\begin{eqnarray*}
x_{1}&=&\frac{qr-\beta_{1}r}{q+1} ,\qquad x_{2}=\frac{q(1-2\rho)-2+2\beta
_{1}\rho
}{2\rho(q+1)}, \\
x_{3}&=&\frac1{\rho}-\beta_{1}(1+r),\qquad x_{4}=\frac2{\rho}-\beta
_{1}(2+3r),\qquad x_{5}=\beta_{1}-\beta_{j+1}.
\end{eqnarray*}

(b) Now, for proving~(\ref{D30}), it clearly follows from~(\ref{D119})
that it
suffices to show that one can choose $q$ and $r$ in such a way that $x_{k}
\geq\beta_{1}/2$ for $k=1,2,3,4$. When $q\rightarrow\infty$ we see that
$x(1)\to x^{\prime}(1)=r$ and $x(2)\to x^{\prime}(2)=\frac{1-2\rho
}{2\rho}$,
so it remains to show that one can choose $r\in(0,1)$ such that
$x^{\prime
}(k)\geq\beta_{1}/2$ for $k=1,2$ and $x_{k}\geq\beta_{1}/2$ for $k=3,4$.
Letting $r$ be bigger than but as close as possible to $\beta_{1}/2$, we
deduce from~(\ref{D31}) that such a choice or $r$ is possible, and the proof
is complete.
\end{pf*}

\section{\texorpdfstring{Proof of Theorem \protect\ref{TD1}}{Proof of Theorem 2}}

(1) In addition to the strengthened Assumption~\ref{A2}, we assume~(\ref{118})
for some $\varepsilon>0$. Theorem~\ref{TD1} says something about the
estimators of $\beta_{i}$ and $A^{i}_{T}$ only when $\beta_{i}>\frac
{\beta
_{1}}2$. Moreover, if~(\ref{31}) holds for the sequence $\beta_{1}%
,\ldots,\beta_{j+1}$, it also holds for the sequence $\beta^{\prime}%
_{1},\ldots,\beta^{\prime}_{j^{\prime}+1}$, where $j^{\prime}=j$ if
$\beta_{j+1}\geq\frac{\beta_{1}}2$ and $j^{\prime}=\sup(i\dvtx\beta_{i}%
>\frac{\beta_{1}}2)$ otherwise, and where $\beta^{\prime}_{i}=\beta
_{i}$ when
$i\leq j^{\prime}$ and $\beta^{\prime}_{j^{\prime}+1}=\beta_{j^{\prime}+1}
\vee\frac{\beta_{1}}2$. Henceforth, upon discarding the indices such that
$\beta_{i}\leq\frac{\beta_{1}}2$, we can assume without loss of generality
that
%
\begin{equation}
\label{127}
\beta_{1}>\cdots>\beta_{j}>\beta_{j+1}=\frac{\beta_{1}}2.
\end{equation}

Under this additional assumption, we have $\beta_{1}-\beta_{j}<1$, and
(\ref{S10}) yields
%
\begin{equation} \label{S1001}
1\leq i\leq k<j \quad\Rightarrow\quad u_{n,i}^{\beta_{i}-\beta_{i+1}} \log\frac
{1}{u_{n,i}} = o( u_{n,k}^{\beta_{i}-\beta_{k+1}}).\vadjust{\goodbreak}
\end{equation}
Moreover, combining~(\ref{31}),~(\ref{KP}) and~(\ref{D30}), we deduce that
%
\begin{equation} \label{S9}%
U(v_{n},\Delta_{n})_{T} = \sum_{i=1}^{j}\frac{A_{T}^{i}}{v_{n}^{\beta
_{i}}%
}+O_{P} (v_{n}^{-\beta/2})
\end{equation}
for any sequence $v_{n}$ such that $v_{n}\leq u_{n}$, and in particular for
the sequences $v_{n}=u_{n,i}$. All of the proof will rely on this, and below
$H_{i}$ is always given by~(\ref{S12}).

(2) We first consider the case $i=1$ when $j>1$. A simple calculation,
based on
(\ref{S9}) applied with $v_{n}=u_{n}$ and $v_{n}=\gamma u_{n}$, yields
that in
restriction to the set $\Omega_{T}$,
\[
\log\bigl(U(v_{n},\Delta_{n})_{T}/U(\gamma v_{n},\Delta_{n})_{T}\bigr) = (
\beta_{1}-H_{1} u_{n}^{\beta_{1}-\beta_{2}}) \log\gamma+o_{P}%
(u_{n}^{\beta_{1}-\beta_{2}}).
\]
This gives the first part of~(\ref{D5}). It also implies that
\begin{eqnarray*}
u_{n}^{\widetilde{\beta}_{1}^{n}}&=&u_{n}^{\beta_{1}}e^{-(\widetilde{\beta}_{1}
^{n}-\beta_{1})\log(1/u_{n})}\\
&=&u_{n}^{\beta_{1}}\bigl( 1+H_{1}u_{n}
^{\beta_{1}-\beta_{2}} \log(1/u_{n})+o_{P}\bigl(u_{n}^{\beta_{1}-\beta_{2} }%
\log(1/u_{n})\bigr)\bigr) .
\end{eqnarray*}
This and~(\ref{S9}) yield the second part of~(\ref{D5}).

(3) Now we suppose that~(\ref{D5}) holds for all $i\leq k-1$, for some
$k\in\{2,\ldots,j-1\}$. We observe that we have the following
identities, for
all $y=(y_{1},\ldots,y_{k+1})$ and $r=1,\ldots,k+1$:
\begin{eqnarray*}
&&
\sum_{l=0}^{k-1}(-1)^{l}\gamma^{-ly_{r}}\sum_{J\in I(k-1,l)} \gamma
^{\sum_{j\in J}y_{j}}\\
&&\qquad=\prod_{l=1}^{k-1}( 1-\gamma^{y_{i}-y_{r}})
=\cases{
0, &\quad if $r\leq k-1$,\cr
G(k,y,\gamma), &\quad if $r=k$,\cr
G^{\prime}(k,y,\gamma), &\quad if $r=k+1$,}
\end{eqnarray*}
where $G(k,y,\gamma)=\prod_{i=1}^{k-1}( 1-\gamma^{y_{i}-y_{k}}) $
and $G^{\prime}(k,y,\gamma)=\prod_{i=1}^{k-1}( 1-\gamma^{y_{i}-y_{k+1}
}) $. Therefore,~(\ref{S9}) applied to $v_{n}=x\gamma^{l}u_{n,k}$ and
the definition of $U^{n}(k,x)$ yield for all $x\geq1$ fixed, and with
$\beta=(\beta_{1},\ldots,\beta_{k+1})$,
\begin{eqnarray*}
U^{n}(k,x) & = & \sum_{r=1}^{k-1}\frac{A_{T}^{r}}{x^{\beta_{r}} u_{n,k}
^{\beta_{r}}} \sum_{l=0}^{k-1}(-1)^{l}( \gamma^{-l\beta_{r}}
-\gamma^{-l\widetilde{\beta}_{r}^{n}}) \sum_{J\in I(k-1,l)}%
\gamma^{ \sum_{j\in J}\widehat{\beta}_{j}^{n}}\\
&&{} +\sum_{r=k}^{k+1}\frac{A_{T}^{r}}{x^{\beta_{r}} u_{n,k}^{\beta_{r}} }%
\sum_{l=0}^{k-1}(-1)^{l}\gamma^{-l\beta_{r}}\sum_{J\in I(k-1,l)}(
\gamma^{\sum_{j\in J}\widetilde{\beta}_{j}^{n}}-\gamma^{ \sum_{j\in J}
\beta_{j}}) \\
&&{} +\frac{A_{T}^{k}}{x^{\beta_{k}} u_{n,k}^{\beta_{k}}} G(k,\beta
,\gamma)+\frac{A_{T}^{k+1}}{x^{\beta_{k+1}} u_{n,k}^{\beta_{k+1}} }%
G^{\prime}(k,\beta,\gamma)+o_{P}(u_{n,k}^{-\beta_{k+1}}).
\end{eqnarray*}
The functions $z\mapsto\gamma^{-lz}$ are $C^{\infty}$. The induction
hypothesis gives $\widetilde{\beta}_{i}^{n}-\beta
_{i}=O_{P}(u_{n,i}^{\beta
_{i}-\beta_{i+1}})$ for $i=1,\ldots,k-1$. Then~(\ref{S1001}) and $\beta
_{i}-\beta_{i+1}>\varepsilon$ allow us to deduce
\begin{eqnarray*}
&\displaystyle 0\leq l\leq k-1, J\in I(k-1,l) \quad\Rightarrow\quad \gamma^{ \sum_{j\in J}
\widetilde{\beta}_{j}^{n}}-\gamma^{ \sum_{j\in J}\beta_{j}} = o_{P}
(u_{n,k}^{\beta_{k-1}-\beta_{k+1}}),&\\
&\displaystyle 1\leq r\leq k-1 \quad\Rightarrow \quad\gamma^{-l\beta_{r}} -\gamma
^{-l\widetilde
{\beta}_{r}^{n}} = o_{P}(u_{n,i}^{\beta_{r}-\beta_{r+1}}) = o_{P}(u_{n,k}
^{\beta_{i}-\beta_{k+1}}).&
\end{eqnarray*}
Therefore we finally obtain
%
\begin{eqnarray} \label{S16}%
U^{n}(k,x) &=& \frac{A_{T}^{k} G(k,\beta,\gamma)}{x^{\beta_{k}}
u_{n,k}^{\beta_{k}}}+\frac{A_{T}^{k+1} G^{\prime}(k,\beta,\gamma)}%
{x^{\beta_{k+1}} u_{n,k}^{\beta_{k+1}}}+o_{P}(u_{n,k}^{-\beta_{k+1}
})\nonumber\\[-8pt]\\[-8pt]
&=& \frac{A_{T}^{k} G(k,\beta,\gamma)}{x^{\beta_{k}} u_{n,k}^{\beta_{k}
}%
} \biggl( 1+\frac{H_{k} \log\gamma}{\gamma^{\beta_{k}-\beta_{k+1}}
-1} (xu_{n,k})^{\beta_{k}-\beta_{k+1}}+o_{P}(u_{n,k}^{\beta_{k} -\beta
_{k+1}%
})\biggr) ,\hspace*{-10pt}\nonumber
\end{eqnarray}
where the last equality comes from the definition of $H_{k}$ in~(\ref{S12}).
Then exactly as in Step~\ref{step2}, a simple calculation shows the first half of
(\ref{D5}) for $i=k$.

For the second part of~(\ref{D5}), and as in Step~\ref{step2}, we first deduce
from the
above that
\[
u_{n,k}^{\widetilde{\beta}_{k}^{n}} = u_{n,k}^{\beta_{k}}\bigl( 1+H_{k}
u_{n,k}^{\beta_{k}-\beta_{k-1}} \log(1/u_{n,k})+o_{P}\bigl(u_{n,k}^{\beta
_{k}-\beta_{k-1}} \log(1/u_{n,k})\bigr)\bigr) .
\]
Therefore it is enough to show that
\[
u_{n,k}^{\beta_{k}}\Biggl( U(u_{n,k})_{T}-\sum_{i=1}^{k-1}\widetilde{\Gamma
}_{i}^{n} u_{n,k}^{-\widetilde{\beta}_{i}^{n}}\Biggr) = A_{T}^{k}
+o_{P}\bigl(u_{n,k}^{\beta_{k}-\beta_{k-1}} \log(1/u_{n,k})\bigr).
\]
In view of~(\ref{S9}) with $v_{n}=u_{n,k}$ this amounts to proving for
$i=1,\ldots,k-1$,
%
\begin{equation} \label{S17}%
\widetilde{\Gamma}_{i}^{n} u_{n,k}^{\beta_{k}-\widetilde{\beta}_{i}^{n}}
-A_{T}^{i} u_{n,k}^{\beta_{k}-\beta_{i}} = o_{P}\bigl(u_{n,k}^{\beta_{k}
-\beta_{k-1}} \log(1/u_{n,k})\bigr).
\end{equation}
The induction hypothesis yields that
\begin{eqnarray*}
u_{n,k}^{\beta_{k}-\widetilde{\beta}_{i}^{n}} &=& u_{n,k}^{\beta_{k}-\beta
_{i}%
}\bigl( 1+O_{P}\bigl(u_{n,i}^{\beta_{i}-\beta_{i+1}} \log(1/u_{n,k})\bigr)\bigr) ,
\\
\widetilde{\Gamma}_{i}^{n} &=& A_{T}^{i}+O_{P}\bigl(u_{n,i}^{\beta_{i}-\beta_{i+1}
} \log(1/u_{n,i})\bigr).
\end{eqnarray*}
Then~(\ref{D5}) readily follows from~(\ref{S10}).

(4) It remains to prove that the variables in~(\ref{D6}) are tight. The
difference with the previous case is that~(\ref{S16}) no longer holds when
$i=j=1$ or $k=j>1$, but it can be replaced by
\[
U^{n}(j,x) = \frac{A_{T}^{j} G(j,\beta,\gamma)}{x^{\beta_{j}} u_{n,j}%
^{\beta_{j}}} \bigl( 1+O_{P}(u_{n,j}^{\beta_{j}-\beta/2})\bigr) .
\]
The rest of the proof goes unchanged [note that $\eta$ in~(\ref{D6}) is
$\beta_{j}-\beta/2$ here].

\section{\texorpdfstring{Proof of Theorem \protect\ref{TD2}}{Proof of Theorem 3}}

We use simplifying notation: a point in $D$ is $\theta=(x_{i},\gamma
_{i})_{1\leq i\leq j}$, and we define the functions $F_{n,l}(\theta
)=\sum_{i=1}^{j}\gamma_{i}/(v_{l}u_{n})^{x_{i}}$. The ``true
value'' of the parameter is $\theta_{0}=(\beta
_{i},\Gamma
_{i})_{1\leq i\leq j}$, the preliminary estimators are $\widetilde
{\theta}
_{n}=(\widetilde{\beta}_{i}^{n},\widetilde{\Gamma}_{i}^{n})_{1\leq
i\leq j}$,
and the final estimators are $\overline{\theta}_{n}=(\overline{\beta}{}^{n}_{i}
,\overline{\Gamma}{}^{n}_{i})_{1\leq i\leq j}$. We set $h_{n} =\log
(1/u_{n})$, and as in the previous proof we can assume~(\ref{127}).

(1) We introduce some specific notation. For $m\geq2$ we set $G_{m}%
=(1,\infty)^{m-1}$, a~point in $G_{m}$ being denoted as $\overline{v}%
=(v_{2},\ldots,v_{m})$. For $1\leq k\leq j$ and $\overline{v}\in
G_{2k}$, and
with the convention $v_{1}=1$, we let $\Sigma(\overline{v})$ be the
$2k\times2k$ matrix with entries
%
\begin{equation}
\label{D2-25}\Sigma(\overline{v})_{l,i} = \cases{
v_{l}^{-\beta_{i}}, &\quad if $1\leq i\leq k$,\cr
v_{l}^{-\beta_{i-k}} \log v_{l}, &\quad if $k+1\leq i\leq2k$.}
\end{equation}
The aim of this step is to show that the set $Z_{k}$ of all $\overline
{v}\in
G_{2k}$ for which the matrix $\Sigma(\overline{v})$ is invertible satisfies
$\lambda_{2k}((Z_{k})^{c})=0$, where $\lambda_{r}$ is the Lebesgue
measure on~$G_{r}$.

When $1\leq m\leq2k$ and $\overline{v}\in G_{2k}$, we denote by
$\mathcal{M}%
_{m}(\overline{v})$ the family of all $m\times m$ sub-matrices of the
$m\times2k$ matrix $(\Sigma( \overline{v})_{l,r}\dvtx 1\leq l\leq m,1\leq
r\leq
2k)$. A key fact is that $\mathcal{M}_{m}(\overline{v})=\mathcal{M}%
_{m}(\overline{v}_{m})$ only depends on the restriction $\overline{v}%
_{m}=(v_{2},\ldots, v_{m})$ of $\overline{v}$ to its first $m-1$ coordinates.
Moreover, $\Sigma(\overline{v})_{1i}$ equals $1$ if $i\leq k$ and $0$
otherwise: so the entries of the first column of any $M\in\mathcal{M}%
_{m}(\overline{v })$ are $0$ or $1$, and $\mathcal{M}^{\prime}_{m}%
(\overline{v})$ denotes the subset of all $M\in\mathcal{M}_{m}(\overline{v})$
for which $M_{1,i}=1 $ for at least one value of $i$. Finally, $H_{m}$ stands
for the set of all $\overline{v}_{m}\in G_{m}$ such that all $M\in
\mathcal{M
}^{\prime}_{m}(\overline{v}_{m})$ are invertible. Since $\mathcal{M}
^{\prime
}_{2k}(\overline{v})$ is the singleton $\{\Sigma(\overline{ v})\}$, we have
$Z_{k}=H_{2k}$.

If $m\geq2$ and $\overline{v}_{m}=(v_{2},\ldots,v_{m})\in G_{m}$ and
$M\in\mathcal{M}_{m}^{\prime}(\overline{v}_{m})$, by expanding along
the last
column, we see that
%
\begin{equation}
\label{B4-1201}
\det(M) = \sum_{i=1}^{k}v_{m}^{\beta_{i}}(a_{i}+a_{k+i}
\log v_{m}),
\end{equation}
where each $a_{r}$ is of the form: either (i) $a_{r}$ is plus or minus
$\det(M_{r})$ for some $M_{r}\in\mathcal{M}_{m-1}(\overline{v}_{m})$
(for $m$
values of $r$) or (ii) $a_{r}=0$ (for the other $2k-m$ values of~$r$). Note
that we can also have $a_{r}=0$ in case (i), and since $M\in\mathcal{M}%
^{\prime}_{m}(\overline{v}_{m})$ there is at least one $a_{r}$ of type (i)
with $M_{r}\in\mathcal{M}^{\prime}_{m-1}(\overline{v}_{m})$.

When at least one $a_{r}$ in~(\ref{B4-1201}) is not $0$, the right-hand
side of
this expression, as a function of $v_{m}$, has finitely many roots only,
because all $\beta_{i}$'s are distinct. Observing that $\mathcal{M}
^{\prime
}_{1}(\overline{v})$ is the $1\times1$ matrix equal to $1$, it follows that,
with $(\overline{v}_{m-1},v_{m})= (v_{2},\ldots,v_{m-1},v_{m})$ when
$\overline{v}_{m-1} =(v_{2},\ldots,v_{m-1})$, and recalling that with our
standing notation $\lambda_{2}$ is the Lebesgue measure on $(1,\infty)$,
%
\begin{eqnarray}
\label{B4-12}\qquad
m=2 &\quad\Rightarrow\quad& \lambda_{2}((H_{2})^{c}) =
0,\nonumber\\[-8pt]\\[-8pt]
m\geq3, \qquad\overline{v}_{m-1}\in H_{m-1} &\quad\Rightarrow\quad& \lambda_{2}\bigl( v_{m}%
\dvtx(\overline{va}_{m-1},v_{m})\notin H_{m}\bigr) = 0.
\nonumber
\end{eqnarray}
Since
\[
\lambda_{m}((H_{m})^{c}) = \int_{G_{m-1}} \lambda_{2}\bigl(v_{m}\dvtx (
\overline{v}_{m-1},v_{m})\notin H_{m}\bigr) \lambda_{m-1}(d \overline
{v}_{m-1}),
\]
which equals $\int_{H_{m-1}} \lambda_{1}(v_{m}\dvtx (\overline{v}
_{m-1},v_{m})\notin H_{m}) \lambda_{m-1} (d\overline{v}_{m-1})$
if
$\lambda_{m-1}((H_{m-1})^{c})=0$, when $m\geq3$, we deduce from~(\ref{B4-12}),
by induction on~$m$, that indeed $\lambda_{m}((H_{m})^{c})=0$ for all
$m=2,\ldots,2k$. Recalling $Z_{k}=H_{2k}$, the result follows.

Since the claim of the theorem holds for all $(v_{2},\ldots,v_{L})$
outside a
$\lambda_{L}$-null set only, and $L\geq2k$, we thus can and will assume below
that the numbers~$v_{l}$ are such that $\overline{v} _{2k}=(v_{2}%
,\ldots,v_{2k}) \in Z_{k}$, hence $\Sigma(\overline{ v}_{2k})$ is invertible,
for all $k=1,\ldots,j$.

(2) Our assumptions on the preliminary estimators yield that the set
$\Omega_{n}$ on which $\Vert\widetilde{\theta}_{i}^{n}-\theta_{0}\Vert
\leq1/u_{n}^{\eta}$ satisfies $\mathbb{P}(\Omega_{n})\rightarrow1$. So below
we argue on the set~$\Omega_{n}$, or equivalently (and more
conveniently) we
suppose $\Omega_{n}=\Omega$. Then~$\overline{\theta}_{n}$ converges pointwise
to $\theta_{0}$, which belongs to all the sets $D_{n}$. Set
\[
y_{i}^{n} = A_{T}^{i}(\overline{\beta}{}^{n}_{i}-\beta_{i}),\qquad z_{i}
^{n} = \overline{\Gamma}{}^{n}_{i}-A_{T}^{i}+y_{i}^{n}h_{n},\qquad a_{i}%
^{n} = |y_{i}^{n}|h_{n}+|z_{i}^{n}|.
\]
We have $a_{i}^{n}\leq2u_{n}^{-\eta}h_{n}$ because $\Omega_{n}=\Omega$. Then
an expansion of $(x_{i},\gamma_{i})\mapsto\gamma
_{i}/(v_{l}u_{n})^{x_{i} }$
around $(\beta_{i},A_{T}^{i})$ yields for all $l$,
%
\begin{equation} \label{D2-5}%
\frac{\overline{\Gamma}_{i}}{(v_{l}u_{n})^{\overline{\beta}_{i}}}-\frac
{A_{T}^{i}}{(v_{l}u_{n})^{\beta_{i}}}= \frac{1}{(v_{l}u_{n})^{\beta_{i}}
}( z_{i}^{n}-y_{i}^{n}\log v_{l}+x_{i,l}^{n}),
\end{equation}
where
\[
|x_{i,l}^{n}| \leq K|y_{i}^{n}|h_{n}(|z_{i}^{n}|+|y_{i}^{n} |) \leq
K|y_{i}^{n}|h_{n}a_{i}^{n} \leq K(a_{i}^{n})^{2}.
\]

Combining~(\ref{31}),~(\ref{KP}) and~(\ref{D30}), we see that
\[
U(v_{l}u_{n},\Delta_{n})_{T}-F_{n,l}(\theta_{0}) = O_{P}(u_{n}^{-\beta_{1}
/2}).
\]
Since $\Phi_{n}(\theta)=\sum_{l=1}^{L}w_{l}( U(v_{l}u_{n},\Delta
_{n})_{T}-F_{n,l}(\theta)) ^{2}$, we deduce
\[
\Phi_{n}(\theta_{0}
)=O_{P}(u_{n}^{-\beta_{1}}).
\]
Since $\theta_{0}\in D_{n}$ and $\overline
{\theta}_{n}$ minimizes $\Phi_{n}$ over $D_{n}$, we also have $\Phi
_{n}(\overline{\theta}_{n})=O_{P}(u_{n}^{-\beta_{1}})$, hence $F_{n,l}
(\theta_{0})-F_{n,j}(\overline{\theta}_{n})=O_{P}(u_{n}^{-\beta
_{1}/2})$ for
all $l$. Using~(\ref{D2-5}), this can be rewritten as
%
\begin{equation}\label{D2-6}%
\sum_{i=1}^{j}\frac{1}{(v_{p}u_{n})^{\beta_{i}}} ( z_{i}^{n}-y_{i}
^{n}\log v_{l}+x_{i,l}^{n}) = O_{P}(u_{n}^{-\beta_{1}/2}).
\end{equation}

(3) Taking $k$ between $1$ and $j$, we consider the $2k$-dimensional vectors
$\zeta(k,n)$ and $\xi(n)$ with components (for $l=1,\ldots,2k$),
\begin{eqnarray*}
\zeta(k,n)_{l} &=& \sum_{i=1}^{k}\frac{1}{(v_{p}u_{n})^{\beta_{i}}} (
z_{i}^{n}-y_{i}^{n}\log v_{l}) ,\vadjust{\goodbreak}\\
\xi(k,n)_{i}&=&\cases{
z_{i}^{n} u_{n}^{-\beta_{i}}, &\quad if $1\leq i\leq k$,\cr
-y_{i-k}^{n} u_{n}^{-\beta_{i-k}}, &\quad if $k+1\leq i\leq2k$.}
\end{eqnarray*}
With matrix notation, and~(\ref{D2-25}), we have $\zeta(k,n)= \Sigma
(\overline{v}_{2k})\xi(k,n)$, hence
%
\begin{equation} \label{D2-7}%
\xi(k,n) = \Sigma(\overline{v}_{2k})^{-1} \zeta(k,n).
\end{equation}

Next, we have
\[
\frac{1}{(v_{l}u_{n})^{\beta_{i}}} \bigl|z_{i}^{n}+\bigl(v^{ \prime}_{g}(\beta
_{i})+v_{g}(\beta_{i})\log\delta_{l}\bigr)y_{i} ^{n}+x_{i,l}^{n}\bigr| \leq
\frac{Ka_{i}^{n}}{u_{n}^{\beta_{i}}},\qquad \frac{|x_{i,l}^{n}|}{(v_{l}%
u_{n})^{\beta_{i}}}\leq\frac{K(a_{i}^{n})^{2}}{u_{n}^{\beta_{i}}},
\]
and hence~(\ref{D2-6}) and $a_{i}^{n}\leq Ku_{n}^{\eta}h_{n}\leq
K/h_{n}^{2} \leq
K$ yield
\[
|\zeta(k,n)_{l}| \leq K\Biggl( \sum_{i=1}^{k-1}(a_{i}^{n})^{2} u_{n}
^{-\beta_{i}}+\frac{a^{n}_{k}}{h_{n}^{2}} u_{n}^{-\beta_{k}}+\sum
_{i=k+1}%
^{j}a_{i}^{n} u_{n}^{-\beta_{i}}\Biggr) +O_{P}(u_{n}%
^{-\beta_{j+1}}).
\]
By~(\ref{D2-7}) the variables $\xi(k,n)_{l}$ satisfy the same estimate. Since
$a_{k}^{n}\leq(|\xi(k,n)_{k}|+ |\xi(k,n)_{2k} |h_{n})u_{n}^{\beta_{k}}$,
\[
a_{k}^{n} \leq Ch_{n}\Biggl( \sum_{i=1}^{k-1}(a_{i}^{n})^{2} u_{n}^{\beta
_{k}-\beta_{i}}+\frac{a^{n}_{k}}{h_{n}^{2}} +\sum_{i=k+1}^{j}a_{i}^{n}
u_{n}^{\beta_{k}-\beta_{i}}\Biggr) +O_{P}(h_{n} u_{n}^{
\beta_{k}-\beta_{j+1}})
\]
for some constant $C$. When $n$ is large enough, $C/h_{n}\leq\frac12$,
and we
deduce
%
\begin{equation}
\label{B4-17}\qquad
a_{k}^{n} \leq2Ch_{n}\Biggl( \sum_{i=1}^{k-1}(a_{i}^{n})^{2}
u_{n}^{\beta_{k}-\beta_{i}} +\sum_{i=k+1}^{j}a_{i}^{n} u_{n}^{\beta
_{k}-\beta_{i}}\Biggr) +O_{P}(h_{n} u_{n}^{\beta_{k}-\beta
_{j+1}}).
\end{equation}

(4) In view of the definition of $y^{n}_{i}$ and $z^{n}_{i}$, to get the
result, and recalling that we assume $\beta_{j+1}=\beta_{1}/2$, it is clearly
enough to prove the existence of a~number $\nu>0$ such that, for all
$i=1,\ldots,j$, we have
%
\begin{equation}
\label{B4-21}
a^{n}_{i} = O_{P}(h_{n}^{\nu} u_{n}^{\beta_{i}-\beta_{j+1}}).
\end{equation}

To this aim, we introduce the following property, denoted
($P_{m,q,r}$), where~$r$ runs through $\{1,\ldots,j\}$ and $m,q\geq1$, and where we use the
notation $\zeta_{r}=\beta_{r}-\beta_{r+1}$:
%
\begin{equation}
\label{B4-22}
i=1,\ldots,r \quad\Rightarrow\quad a_{i}^{n} = O_{P}\bigl( h_{n}%
^{m} (u_{n}^{\beta_{i}-\beta_{r}+q\zeta_{r}} +u_{n}^{\beta_{i}-\beta_{r+1}})\bigr).
\end{equation}
Since $a_{i}^{n}\leq K$, applying~(\ref{B4-17}) with $k=1$ yields $a_{1}
^{n}=O_{P}(h_{n}u_{n}^{\beta_{1}-\beta_{2}})$, which is ($P_{1,1,1}$).

Next, we suppose that ($P_{m,q,r}$) holds for some $r<j$, and for some
$m,q\geq1$. Letting first $k=r+1$, we deduce from~(\ref{B4-17}) that, since
again $a^{n}_{i}\leq K$,
%
\begin{eqnarray}
\label{B4-23}%
a^{n}_{k} & = & O_{P}\Biggl(h_{n}^{1+2m}\sum_{i=1}^{k-1}\bigl( u_{n}^{ \beta
_{k}-\beta_{i}+2(\beta_{i}-\beta_{r}+q\zeta_{r})} +u_{n}^{\beta_{i}%
-\beta_{r+1}} \bigr)\nonumber\\
& &\hspace*{70pt}{} +h_{n}\sum_{i=k+1}^{j} u_{n}^{\beta_{k}-\beta_{i}}+h_{n}%
u_{n}^{\beta_{k}-\beta_{j+1}}\Biggr)\\
& = & O_{P}\bigl(h_{n}^{1+2m} ( u_{n}^{\beta_{k}-\beta
_{r}+2q\zeta_{r}} +u_{n}^{\zeta_{r}}+u_{n}^{\beta_{k}-\beta_{r+2}})\bigr),
\nonumber
\end{eqnarray}
where the last line holds because $k=r+1$ and $h_{n}>1$ for $n$ large enough
and the sequence $\beta_{i}$ is decreasing. This in turn implies, for $k=r+1$
again,
%
\begin{equation}
\label{B4-2333}
a^{n}_{k}=O_{P}\bigl(h_{n}^{r+2-k+2m} ( u_{n}^{\beta
_{k}-\beta_{r}+2q \zeta_{r}}+u_{n}^{\beta_{k}-\beta_{r+1}})\bigr).
\end{equation}

Then, exactly as above, we apply~(\ref{B4-17}) with $k=r$, and~(\ref{B4-22})
and also~(\ref{B4-2333}) with $k=r+1$, to get that~(\ref{B4-2333})
holds for
$k=r$ as well. Repeating the argument, a~downward induction yields that indeed
(\ref{B4-2333}) holds for all $k$ between $1$ and $r+1$. Thus~(\ref{B4-22})
holds with $q$ and $m$ substituted with $2q$ and $r+1+2m$. Hence ($P_{m,q,r}$)
implies ($P_{r+1+2m,2q,r}$). Since obviously ($P_{m,q,r}$) $\Rightarrow$
($P_{m,q^{\prime},r}$) for any $q^{\prime}\in[1,q]$, by a repeated use
of the
previous argument we deduce that if ($P_{m,1,r}$) holds for some $m\geq1$,
then for any $q^{\prime}\geq1$ we can find $m(q^{\prime})\geq1$ such that
($P_{m(q^{\prime}),q^{\prime},r}$) holds as well.

Now, assuming ($P_{m,q,r}$) for some $m,q,r$, we take $q^{\prime}=\frac
{\zeta_{r+1}}{2\zeta_{r}} \vee1$ and $m^{\prime}=m(q^{\prime})$. What precedes
yields ($P_{m^{\prime},q^{\prime},r}$), hence~(\ref{B4-23}) holds for all
$k\leq r+1$, with $q^{\prime}$ and $m^{\prime}$. In view of our choice of
$q^{\prime}$, this implies that ($P_{r+1+m^{\prime},1,r+1}$) holds. Since
($P_{1,1,1}$) holds, we see by induction that for any $r\leq j$ there exists
$m_{r}\geq1$ such that ($P_{m_{r},1,r}$) holds.

It remains to apply~(\ref{B4-22}) with $r=j$ and $m=m_{r}$ and $q=1$,
and we
get~(\ref{B4-21}) with $\nu=m_{j}$. This completes the proof.

\section{\texorpdfstring{Proof of Theorem \lowercase{\protect\ref{theo-fisher}}}{Proof of Theorem 4}}

The proof of Theorem~\ref{theo-fisher} is contained in the supplemental
article [\citet{yacjacod11asupp}].%
\end{appendix}

\begin{supplement}
\stitle{Supplement to ``Identifying the successive Blumenthal--Getoor
indices of a discretely observed process''}
\slink[doi]{10.1214/12-AOS976SUPP} 
\sdatatype{.pdf}
\sfilename{aos976\_supp.pdf}
\sdescription{This supplement contains the proof of Theorem~\ref{theo-fisher}.}
\end{supplement}

%

\printaddresses

\end{document}